\documentclass[11pt]{article}

\usepackage[T2A]{fontenc}
\usepackage[cp1251]{inputenc}
\usepackage{amsthm}
\usepackage{amsfonts}
\usepackage{eufrak}
\usepackage{amssymb}
\usepackage{amsmath}
\usepackage{cite}


\begin{document}
\newtheoremstyle{mytheorem}
  {\topsep}   
  {\topsep}   
  {\itshape}  
  {}       
  {\bfseries} 
  {  }         
  {5pt plus 1pt minus 1pt} 
  { }          
\newtheoremstyle{myremark}
  {\topsep}   
  {\topsep}   
  {\upshape}  
  {}       
  {\bfseries} 
  {   }         
  {5pt plus 1pt minus 1pt} 
  { }          
\theoremstyle{mytheorem}
\newtheorem{theorem}{Theorem}[section]
 \newtheorem{theorema}{Theorem}
  \newtheorem*{heyde*}{Heyde's theorem}
 \newtheorem*{a*}{The Skitovich--Darmois theorem}
 \newtheorem*{b*}{The Skitovich--Darmois theorem for the group $\mathbb{R}\times \mathbb{Z}(2)$}
\newtheorem{proposition}{Proposition}[section]
\newtheorem{lemma}{Lemma}[section]
\newtheorem{corollary}{Corollary}[section]
\newtheorem{definition}{Definition}[section]
\theoremstyle{myremark}
\newtheorem{remark}{Remark}[section]
\bigskip
\centerline{\textbf{On Heyde's  theorem for locally compact }}
\centerline{\textbf{ Abelian groups containing    elements  of order 2}}


\bigskip

\centerline{\textbf{G.M. Feldman}}


\bigskip

\noindent{\bf Abstract}

\medskip

\noindent According to the well-known Heyde theorem the class of Gaussian distributions on the real line  is characterized by the symmetry of the conditional distribution of one linear form
 of independent random variables given  the other.
 We study    analogues of this theorem for some locally compact Abelian groups $X$   containing  an   element of order 2. We prove that if $X$ contains an element of order 2,  this  can lead to the fact that   a wide class of non-Gaussian distributions on $X$ is characterized by the symmetry of the conditional distribution of one linear form  given the other. In so doing coefficients of linear forms are topological automorphisms of $X$.

\bigskip


\bigskip


\section{  Introduction}

In   \cite{He}, see also  \cite[\S 13.4.1]{KaLiRa},  C.C. Heyde proved the following characterization  theorem for the class of Gaussian distributions  on the real line.
\begin{heyde*} \label{30}    Let $\xi_j$, $j = 1, 2,\dots, n$, $n \ge 2,$ be independent random variables with distributions  $\mu_j$. Let $a_j$, $b_j$ be nonzero real numbers satisfying the conditions
 $b_ia_i^{-1} + b_ja_j^{-1}\ne 0$   for all $i, j$.
 Assume that the conditional distribution of the linear form
$L_2 = b_1\xi_1 + \dots + b_n\xi_n$ given $L_1 = a_1\xi_1 + \dots + a_n\xi_n$ is symmetric.
Then  all distributions $\mu_{j}$   are Gaussian, possibly degenerate.
\end{heyde*}

 In the case when independent random variables
 take values in a locally compact Abelian group   $X$,
and coefficients of   linear forms are topological automorphisms of
 $X$, some analogues of Heyde's theorem  were studied in  \cite{Fe2, Fe4, Fe3, Fe20bb, Fe6,      F, FeTVP1, FeTVP, F_solenoid, My2, My1, MiFe1}, see also   \cite[Chapter VI]{Fe5}. In this article we continue  to study group analogues of    Heyde's   theorem. We consider  locally compact Abelian groups that contain  an  element of order 2. Elements of order 2 play an important  role in   Heyde's theorem.  We prove that  if a group $X$ contains an element of order 2,  it can lead  to the fact that   a wide class of non-Gaussian distributions on $X$ is characterized by the symmetry of the conditional distribution of one linear form of independent random variables given the other. The proofs of the corresponding theorems are reduced to solving some functional equations in the class of continuous positive definite functions on the character group of the    group $X$.

Let $X$ be a second countable locally compact Abelian group. Denote by ${\rm Aut}(X)$ the group
of topological automorphisms of $X$, and by  $I$ the identity automorphism of a group.  Denote by $Y$ the character
group of the group $X$, and by  $(x,y)$ the value of a character $y \in Y$ at an
element $x \in X$.  If $K$ is a closed subgroup of $X$, denote by
 $A(Y, K) = \{y \in Y: (x, y) = 1$ \mbox{ for all } $x \in K \}$
its annihilator. Let $X_1$ and
$X_2$ be locally compact Abelian groups with the character groups
$Y_1$ and $Y_2$ respectively. Let
 $\alpha\colon X_1\to X_2$
be a continuous homomorphism. The adjoint homomorphism $\tilde\alpha\colon Y_2\to Y_1$
is defined by the formula $(\alpha x_1,
y_2)=(x_1 , \tilde\alpha y_2)$ for all $x_1\in X_1$, $y_2\in
Y_2$.    Note that $\alpha\in {\rm Aut}(X)$
if and only if $\tilde\alpha\in {\rm Aut}(Y)$. Put $X_{(2)}=\{x\in X: 2x=0\}$, $X^{(2)}=\{2x: x\in X\}$.
Denote by   $\mathbb{R}$ the group of real numbers, by
 $\mathbb{Z}$ the group of integers,  by ${\mathbb Z}(2)=\{0, 1\}$ the group of  residue  classes modulo $2$, and by $\mathbb{T}$ the circle group, i.e. the multiplicative group of all complex numbers with absolute value 1.

  Denote by ${\rm M}^1(X)$ the
convolution semigroup of probability distributions on the group $X$. Let
${\mu\in {\rm M}^1(X)}$.
Denote by
$$
\hat\mu(y) =
\int_{X}(x, y)d \mu(x), \quad y\in Y,$$  the characteristic function  (Fourier transform) of
the distribution  $\mu$, and by $\sigma(\mu)$ the support of $\mu$. The characteristic function of a signed measure on $X$ is defined in the same way. Define the distribution $\bar \mu \in {\rm M}^1(X)$ by the formula
 $\bar \mu(B) = \mu(-B)$ for any Borel  subset $B$ of $X$.
Then $\hat{\bar{\mu}}(y)=\overline{\hat\mu(y)}$. If $G$ is a Borel subgroup of $X$, we denote by ${\rm M}^1(G)$ the subsemigroup of  ${\rm M}^1(X)$  of distributions concentrated at $G$.
Denote by $E_x$  the degenerate distribution
 concentrated at an element $x\in X$ and by $m_K$ the Haar distribution on a compact subgroup
 $K$ of the group $X$.
The characteristic function of a distribution
$m_K$ is of the form
\begin{equation}\label{08_1}
\hat m_K(y)=
\begin{cases}
1, & \text{\ \ }\   y\in A(Y, K),
\\  0, & \text{\ \ }\ y\not\in
A(Y, K).
\end{cases}
\end{equation}

A distribution  $\gamma\in {\rm M}^1(X)$  is called Gaussian
(\cite[Chapter IV, \S 6]{Pa}, see also \cite{Fe1979})
if its characteristic function is represented in the form
\begin{equation}\label{f1}
\hat\gamma(y)= (x,y)\exp\{-\varphi(y)\}, \quad  y\in Y,
\end{equation}
where $x \in X$, and $\varphi(y)$ is a continuous nonnegative function
on the group $Y$
 satisfying the equation
 \begin{equation}\label{f2}
\varphi(u + v) + \varphi(u
- v) = 2[\varphi(u) + \varphi(v)], \quad u,  v \in
Y.
\end{equation}
Denote by $\Gamma(X)$ the set of Gaussian distributions on
the group $X$. According to the definition the degenerate distributions are Gaussian. Note that the support of a Gaussian distribution is a coset of a connected subgroup of $X$.

\section{   Heyde's characterization theorem for the group $\mathbb{R}\times \mathbb{Z}(2)$}

The simplest locally compact Abelian group containing an element of order 2 where there are nondegenerate Gaussian distributions  is the group $X=\mathbb{R}\times \mathbb{Z}(2)$. In this section we study   an analogue of Heyde's theorem for the group $X$. We   prove that   the class of distributions on $X$  which is characterized by the symmetry of the conditional distribution of one linear form of independent random variables given the other, generally speaking,   is much wider than the class $\Gamma(\mathbb{R})*{\rm M}^1(\mathbb{Z}(2))$.

Denote by $x=(t, k)$, $t\in \mathbb{R}$, $k\in \mathbb{Z}(2)$,    elements of the group $X$.
 The  group $Y$  is topologically isomorphic to the group $X$, i.e. $\mathbb{R}\times \mathbb{Z}(2)$.
Denote by $y=(s, l)$, $s\in \mathbb{R}$,
 $l\in \mathbb{Z}(2)$,
elements of the group $Y$.
Every automorphism  $a\in {\rm Aut}(X)$ is of the form
 $a(t, k)=(c_a t, k)$, where $c_a$ is a nonzero real number.
In order not to complicate the notation, if $a\in {\rm Aut}(X)$, then we will write
$a(t, k)=(a t, k)$  and assume that
$a\in \mathbb{R}$, $a\ne 0$.  Obviously, $a=\tilde a$.

Let $\mu$ be a distribution on the group $X$, and let $\mu\in \Gamma(\mathbb{R})*{\rm M}^1(\mathbb{Z}(2))$, i.e. $\mu=\gamma*\omega$, where $\gamma\in\Gamma(\mathbb{R})$, $\omega\in {\rm M}^1(\mathbb{Z}(2))$.
Then the characteristic function of the distribution  $\mu$ is of the form
$$
\hat\mu(s, l) = \begin{cases}\exp\{-\sigma s^2+i\beta s\},
&s\in \mathbb{R}, \ l=0,\\ \kappa\exp\{-\sigma s^2+i\beta s\},
&s\in \mathbb{R}, \   l=1,
\\
\end{cases}
$$
where $\sigma \ge 0$, $\beta, \kappa\in \mathbb{R}$, $|\kappa|\le 1$.
We consider now a class of distributions on the group
$X$  which is wider than the class $\Gamma(\mathbb{R})*{\rm M}^1(\mathbb{Z}(2))$. This class of  distributions    first appeared in \cite{F_solenoid} in connection with the study of    Heyde's theorem on $\text{\boldmath $a$}$-adic solenoids.
We need the following statement.
\begin{lemma}  [\!\!{\protect\cite{F_solenoid}}] \label{lem6} Let   $X=\mathbb{R}\times \mathbb{Z}(2)$, and let $f(s, l)$ be a function on the group $Y$
of the form
\begin{equation}\label{y6}
f(s, l) = \begin{cases}\exp\{-\sigma s^2+i\beta s\}, &s\in \mathbb{R}, \ l=0,\\ \kappa\exp\{-\sigma' s^2+i\beta's\}, &s\in \mathbb{R}, \   l=1,
\\
\end{cases}
\end{equation}
where $\sigma\ge 0$,  $\sigma'\ge 0$,  $\beta, \beta', \kappa\in \mathbb{R}$. Then
$f(s, l)$ is the characteristic function of a signed measure  $\mu$ on the group $X$. Moreover,
   $\mu$ is a probability distribution if and only if
either  $0<\sigma'<\sigma$ and $0<|\kappa|\le\sqrt{\sigma'\over \sigma}\exp\left\{-{(\beta-\beta')^2\over 4(\sigma-\sigma')}\right\}$
or $\sigma=\sigma'$, $\beta=\beta'$ and  $|\kappa|\le 1$. In the last case $\mu\in \Gamma(\mathbb{R})*{\rm M}^1(\mathbb{Z}(2))$.
\end{lemma}
\begin{definition}\label{d1} {Let $X=\mathbb{R}\times \mathbb{Z}(2)$, and let $\mu\in {\rm M}^1(X)$. We say that    $\mu\in\Theta$  if the characteristic function   $\hat\mu(s, l)$ is represented in the form $(\ref{y6})$, where either  $0<\sigma'<\sigma$ and $0<|\kappa|\le\sqrt{\sigma'\over \sigma}\exp\left\{-{(\beta-\beta')^2\over 4(\sigma-\sigma')}\right\}$
or $\sigma=\sigma'$, $\beta=\beta'$ and  $|\kappa|\le 1$.}
\end{definition}

Introduce into consideration another class     of probability distributions on the group $X$.
Let $\mu\in {\rm M}^1(X)$. Define a distribution   $\mu_{_{\mathbb{R}}}\in {\rm M}^1(\mathbb{R})$
by the formula $\mu_{_{\mathbb{R}}}(E)=\mu(E\times{\mathbb{Z}(2))}$, where $E$ is a Borel subset of   $\mathbb{R}$.
\begin{definition}\label{d2} {Let $X=\mathbb{R}\times \mathbb{Z}(2)$, and let $\mu\in {\rm M}^1(X)$. We say that
$\mu\in\Lambda$ if $\mu_{_{\mathbb{R}}}\in \Gamma(\mathbb{R})$. In other words, $\mu\in\Lambda$ if and only  if $\hat\mu(s, 0)$ is the characteristic function of  a Gaussian distribution on the real line, possibly degenerate.}
\end{definition}
Obviously, all distributions $\mu\in\Lambda$ can be obtained as  follows. Let  $\gamma$  be a Gaussian distribution on $\mathbb{R}$ with the characteristic function
$\hat\gamma(s)=e^{-\sigma s^2+i\beta s}$, $s\in \mathbb{R}$, where $\sigma\ge 0$,    $\beta\in \mathbb{R}$.
Let $\gamma=\lambda+\omega$, where $\lambda$ and $\omega$ are measures on $\mathbb{R}$.
Define a distribution $\mu\in {\rm M}^1(X)$ as follows
$$
\mu(E\times \{k\})=\begin{cases}\lambda(E), \ \   k=0,\\ \omega(E),   \ \  k=1,
\end{cases}
$$
where $E$ is a Borel subset of  $\mathbb{R}$. Then $\mu\in\Lambda$,
and $\mu_{_{\mathbb{R}}}=\gamma$.
We note that the classes $\Theta$ and $\Lambda$ are   subsemigroups of ${\rm M}^1(X)$, and $\Theta\subset\Lambda$.

The following statement can be viewed as an  analogue of Heyde's theorem
for  independent random variables with values
in the group $X=\mathbb{R}\times \mathbb{Z}(2)$.
\begin{theorem}\label{th1}
Let $X=\mathbb{R}\times \mathbb{Z}(2)$, and let $a_j$, $b_j$, $j = 1, 2,\dots, n$, $n \ge 2,$ be topological automorphisms of the group $X$
satisfying the conditions $b_ia_i^{-1} + b_ja_j^{-1}\ne 0$  for all $i, j$. Let
$\xi_j$ be independent random variables with values in     $X$ with distributions $\mu_j$.
Assume that the conditional distribution of the linear form
$L_2 = b_1\xi_1 + \dots + b_n\xi_n$ given $L_1 = a_1\xi_1 + \dots + a_n\xi_n$ is symmetric.
 Then the following alternative holds.

{\rm  1}. All distributions $\mu_{j}$ belong to the class $\Lambda$, and at least one of the distributions  $\mu_j$ is represented in the form
$\mu_{j}=\gamma_j*m_{\mathbb{Z}(2)},$
where $\gamma_j\in\Gamma(\mathbb{R})$.

{\rm 2}.  All distributions $\mu_{j}$ belong to the class
 $\Theta$, and their characteristic functions do not vanish.
\end{theorem}
To prove Theorem \ref{th1}  we need some lemmas.
\begin{lemma} [\!\!{\protect\cite[Lemma 16.1]{Fe5}}] \label{lem1}
Let $X$ be a locally compact Abelian group, and let   $\alpha_j$, $\beta_j$, $j = 1, 2,\dots, n$, $n \ge 2,$ be topological automorphisms of the group $X$.  Let
$\xi_j$ be independent random variables with values in
       $X$  and distributions $\mu_j$.  The conditional distribution of the linear form
$L_2 = \beta_1\xi_1 + \dots + \beta_n\xi_n$ given $L_1 = \alpha_1\xi_1 + \dots + \alpha_n\xi_n$ is symmetric if and only
 if the characteristic functions
 $\hat\mu_j(y)$ satisfy the equation
\begin{equation}\label{42}
\prod_{j=1}^n \hat\mu_j(\tilde \alpha_j u+\tilde \beta_j v )=
\prod_{j=1}^n \hat\mu_j(\tilde \alpha_j u-\tilde \beta_j v), \quad u, v \in Y.
\end{equation}
\end{lemma}
It is convenient for us to formulate as   lemmas the following   well known statements  (see e.g.   \cite[\S 2]{Fe5}).
\begin{lemma}\label{lem2}  Let $X$ be a locally compact Abelian group,   and let $\mu\in{\rm
M}^1(X)$. Then the set $E=\{y\in Y:\
\hat\mu(y)=1\}$ is a closed subgroup of $Y$,  and $\sigma(\mu)\subset
A(X,E)$.
\end{lemma}
\begin{lemma}\label{lem11}  Let $X$ be a locally compact Abelian group, and let  $G$ be a Borel subgroup of $X$. Assume
$\mu\in {\rm M}^1(G)$,  and $\mu=\mu_1*\mu_2$, where $\mu_j\in {\rm
M}^1(X)$. Then the distributions $\mu_j$ may be replaced by their shifts  $\lambda_1=\mu_1*E_x$, $\lambda_2=\mu_2*E_{-x}$, $x\in X$, in such a manner that
$\lambda_j\in {\rm
M}^1(G)$, $j=1, 2$.
\end{lemma}
We   formulate as a lemma the following   statement
which  can be easily verified (see e.g.  \cite[Lemma 6.9]{Fe9}).
\begin{lemma}\label{lem7}  Let $X=\mathbb{R}\times G$, where $G$ is a locally compact Abelian group. Let   $H$ be the character group  of the group  $G$.
Denote by $(s, h)$, $s\in\mathbb{R}$, $h\in H$, elements of the group $Y$. Let
$\mu\in {\rm M}^1(X)$.  Assume that the characteristic function $\hat\mu(s, 0)$ can be extended to the complex plane   as an entire function in $s$. Then for each fixed $h\in H$ the function $\hat\mu(s, h)$ can be also extended to the complex plane    as an entire function in $s$. Moreover, the following inequality
   \begin{equation}\label{8_16}
\max_{|s|\le r}|\hat\mu(s, h)|\le \max_{|s|\le r}|\hat\mu(s, 0)|, \quad h\in H,
\end{equation}
holds.
\end{lemma}
We  need the following particular case of a Linnik--Rao result  (\!\!\cite[Lemma 1.5.1]{KaLiRa}).
\begin{lemma}\label{lem7_1}Let $\varepsilon>0$. Consider the equation
$$
\sum_{j=1}^n\psi_j(s_1+c_js_2)=0,
$$
where $|s_1|<\varepsilon$, $|s_2|<\varepsilon$,  all numbers $c_j$ are pairwise distinct, and the complex valued functions $\psi_j(s)$ in real variable $s$ are continuous. Then the functions $\psi_j(s)$ are polynomials in a neighborhood of zero.
\end{lemma}
We  need the fact that if a distribution with a nonvanishing characteristic function belongs to the class   $\Theta$, then its divisors also belong to the class $\Theta$.
\begin{lemma}  [\!\!{\protect\cite{F_solenoid}}] \label{lem8} Let  $X=\mathbb{R}\times \mathbb{Z}(2)$. Suppose that a distribution $\mu$ on the group $X$ with a nonvanishing characteristic function belongs to the class $\Theta$. If $\mu=\mu_1*\mu_2,$ where $\mu_j\in {\rm M}^1(X)$, then $\mu_j\in \Theta$, $j=1, 2$.
\end{lemma}
{\it Proof of Theorem $\ref{th1}$} \ Passing to the new random variables
 $\zeta_j = a_j \xi_j$ we reduce the proof of the theorem to the case when    $L_1$ and $L_2$
are of the form $L_1 = \xi_1 + \dots + \xi_n$ and $L_2 =b_1 \xi_1 +
\dots + b_n\xi_n$, and topological automorphisms $b_j$ satisfy  the conditions $b_i + b_j\ne 0$ for all $i, j$.
By Lemma \ref{lem1}, the characteristic functions  $\hat\mu_j(s, l)$ satisfy Eq.  (\ref{42})  which takes the form
\begin{equation}\label{05_08_1}
\prod_{j=1}^n \hat\mu_j( s_1+ b_js_2, l_1+l_2)= \prod_{j=1}^n \hat\mu_j( s_1-b_js_2, l_1+l_2), \quad s_1, s_2\in \mathbb{R},
 \ l_1, l_2\in \mathbb{Z}(2).
\end{equation}
Substituting $l_1=l_2=0$ in
(\ref{05_08_1}), we get
\begin{equation}\label{09_08_1}
\prod_{j=1}^n \hat\mu_j( s_1+ b_js_2, 0)= \prod_{j=1}^n \hat\mu_j( s_1-b_js_2, 0), \quad s_1, s_2\in \mathbb{R}.
\end{equation}
Taking into account Lemma \ref{lem1} and Heyde's theorem, (\ref{09_08_1})  implies      that
\begin{equation}\label{05_08_2}
\hat\mu_j(s, 0)=e^{-\sigma_j s^2+i\beta_j s}, \quad s\in \mathbb{R}, \   j = 1, 2,\dots, n,
\end{equation}
where $\sigma_j\ge 0$,  $\beta_j\in \mathbb{R}$. Substituting (\ref{05_08_2}) into  (\ref{09_08_1}), we obtain
\begin{equation}\label{11}
\sum_{j=1}^nb_j\sigma_j=0, \quad \sum_{j=1}^nb_j\beta_j=0.
\end{equation}
It follows from  (\ref{05_08_2}) that   $\mu_j\in \Lambda$, $j = 1, 2,\dots, n$. Denote by $\gamma_j$ the Gaussian distribution on the group  $X$ with the characteristic function
\begin{equation}\label{08_2}
\hat\gamma_j(s, l)=\hat\mu_j(s, 0), \quad s\in \mathbb{R}, \
 l\in \mathbb{Z}(2), \ j = 1, 2,\dots, n.
\end{equation}
Since $\hat\gamma_j(0, l)=1$ for all $l\in \mathbb{Z}(2)$, by Lemma \ref{lem2},  $\sigma(\gamma_j)\subset A(X, \mathbb{Z}(2))=\mathbb{R}$, i.e. $\gamma_j\in \Gamma(\mathbb{R})$. Note that, since $A(Y, \mathbb{Z}(2))=\mathbb{R}$, it follows from  (\ref{08_1}) that
\begin{equation}\label{09_1}
\hat m_{\mathbb{Z}(2)}(s, l)=
\begin{cases}
1, & s\in \mathbb{R}, \   l=0,
\\  0, & s\in \mathbb{R}, \ l=1.
\end{cases}
\end{equation}
Substituting $l_1=1$, $l_2=0$ in (\ref{05_08_1}), we obtain
\begin{equation}\label{08_4}
\prod_{j=1}^n \hat\mu_j( s_1+ b_js_2, 1)= \prod_{j=1}^n \hat\mu_j( s_1-b_js_2, 1), \quad s_1, s_2\in \mathbb{R}.
\end{equation}

Assume that at least for one $j=j_0$ we have  $\hat\mu_{j_0}(0, 1)=0$. Substituting $s_1=b_{j_0}s$, $s_2=s$ in (\ref{08_4}), we receive
\begin{equation}\label{08_5}
\prod_{j=1}^n \hat\mu_j((b_{j_0}+b_j)s , 1)= 0, \quad s\in \mathbb{R}.
\end{equation}
By Lemma \ref{lem7}, it follows from  (\ref{05_08_2})  that   $\hat\mu_j(s, 1)$, $j = 1, 2,\dots, n$, are entire functions. Since by the assumption $b_i + b_j\ne 0$ for all $i, j$, we find from (\ref{08_5}) that $\hat\mu_j(s, 1)=0$ for all $s\in \mathbb{R}$, for at least one of the characteristic functions $\hat\mu_j(s, l)$.
Taking into account (\ref{08_2}) and (\ref{09_1}),   this characteristic function can be represented in the form
$\hat\mu_j(s, l)=\hat\gamma_j(s, l)\hat m_{\mathbb{Z}(2)}(s, l)$. Hence,
$\mu_{j}=\gamma_j*m_{\mathbb{Z}(2)}$.

 Assume now that $\hat\mu_j(0, 1)\ne 0$ for all  $j = 1, 2,\dots, n$.
Then there exists  $\varepsilon>0$ such that $\hat\mu_j(s, 1)\ne 0$, $j = 1, 2,\dots, n$, for $|s|<\varepsilon$. Put
$$\psi_j(s)=\log\hat\mu_j(s, 1), \quad |s|<\varepsilon, \quad  j = 1, 2,\dots, n,
$$
(we consider the principal branch of the logarithm). It follows from  (\ref{08_4}) that in a   neighbourhood of zero the continuous functions $\psi_j(s)$ satisfy the equation
\begin{equation}\label{08_10_1}
\sum_{j=1}^n(\psi_j(s_1+b_js_2)-\psi_j(s_2-b_js_2))=0, \quad s_1, s_2\in \mathbb{R}.
\end{equation}

Assume first that all numbers $b_j$ are pairwise distinct. Since by assumption $b_i + b_j\ne 0$ for all $i, j$, we may apply Lemma \ref{lem7_1}, and it follows from (\ref{08_10_1}) that in a  neighbourhood of zero the continuous functions $\psi_j(s)$  are polynomials. Since   $\hat\mu_j(s, 1)$ are entire functions, we have the representation
$$
\hat\mu_j(s, 1)=e^{\psi_j(s)}, \quad  j = 1, 2,\dots, n,
$$
in the  complex plane.
By Lemma \ref{lem7}, it follows from (\ref{8_16}) and (\ref{05_08_2})  that the entire functions
$\hat\mu_j(s, 1)$ are  of at most order 2. For this reason the degrees of the polynomials $\psi_j(s)$ are at most 2.
Taking into account that $\hat\mu_j(-s, 1)=\overline{\hat\mu_j(s, 1)}$ and $|\hat\mu_j(s, 1)|\le 1$ for all $s\in \mathbb{R}$, we get from (\ref{8_16}) and (\ref{05_08_2}) the representation
\begin{equation}\label{08_7}
\hat\mu_j(s, 1)=\kappa_je^{-\sigma_j' s^2+i\beta_j's}, \quad s\in \mathbb{R}, \   j = 1, 2,\dots, n,
\end{equation}
where $0\le\sigma_j'\le\sigma_j$, $\beta_j', \kappa_j\in \mathbb{R}$, $|\kappa_j|\le
1$, $\kappa_j\ne 0$.
 It follows from (\ref{05_08_2}) and (\ref{08_7})  that $\mu_j\in \Theta$, $j = 1, 2,\dots, n$, and all the characteristic functions $\hat\mu_j(s, l)$ do not vanish. Hence, in the case when all numbers $b_j$ are pairwise distinct, the theorem is proved.

 Get rid of the restriction that all numbers $b_j$ are pairwise distinct. Assume that $b_1=\dots=b_{j_1}$, $b_{j_1+1}=\dots=b_{j_2}$, $\dots$, $b_{j_{l-1}+1}=\dots=b_{j_l}=b_n$, where the numbers  $b_{j_1}$, $b_{j_2}$, $\dots$, $b_{j_l}$ are pairwise distinct. Put $\eta_1=\xi_1+\dots+\xi_{j_1}$, $\eta_2=\xi_{j_1+1}+\dots+\xi_{j_2}$, $\dots$, $\eta_l=\xi_{j_{l-1}+1}+\dots+\xi_{j_l}$. Then $\eta_1$, $\dots$, $\eta_l$ are independent random variables with values in the group  $X$ and distributions $\nu_1=\mu_1*\dots*\mu_{j_1}$, $\nu_2=\mu_{j_1+1}*\dots*\mu_{j_2}$, $\dots$, $\nu_l=\mu_{j_{l-1}+1}*\dots*\mu_{j_l}$.
 The linear forms $L_1$ and $L_2$ of independent random variables $\eta_j$ have the form
 $L_1=\eta_1+\dots+\eta_{l}$ and $L_2=b_{j_1}\eta_1+\dots+b_{j_l}\eta_{l}$. Since the numbers $b_{j_1}$, $b_{j_2}$, $\dots$ ,$b_{j_l}$ are pairwise distinct,   as proved above,
 $\nu_j\in \Theta$, $j = 1, 2,\dots, l$. Obviously, the characteristic functions of the distributions $\nu_j$ do not vanish.  By Lemma \ref{lem8}, this implies  that $\mu_j\in \Theta$, $j = 1, 2,\dots, n$, and the characteristic functions $\hat\mu_j(s, l)$ do not vanish.
 $\blacksquare$
 \begin{corollary}\label{co1}  Assume that in Theorem $\ref{th1}$ $a_j=1$ and   $b_j>0$, $j = 1, 2,\dots, n$.  Then some shifts of  the distributions   $\mu_j$
are supported in $\mathbb{Z}(2)$.
\end{corollary}
{\it Proof} \ It follows from (\ref{11}) that   $\sigma_1=\dots=\sigma_n=0$. Put $\nu_j=\mu_j*\bar \mu_j$. Then $\hat\nu_j(s, 0)=1$ for all $s\in \mathbb{R}$. By Lemma  \ref{lem2}, this implies that $\sigma(\nu_j)\subset A(X, \mathbb{R})=\mathbb{Z}(2)$,     $j=1, 2, \dots, n$.
The   required assertion follows from Lemma  \ref{lem11}. $\blacksquare$
\begin{remark}\label{re1}  Let $X=\mathbb{R}\times \mathbb{Z}(2)$, and let $b_j\in {\rm Aut}(X)$, $j=1, 2\dots, n$. Then we can construct distributions $\mu_j$ such as in   items  1 and    2    of Theorem $\ref{th1}$ with the following property:
 if
 $\xi_j$  are independent random variables with values in the group
       $X$  and distributions $\mu_j$, then the conditional distribution of the linear form
$L_2 = b_1\xi_1 + \dots + b_n\xi_n$ given $L_1 = \xi_1 + \dots + \xi_n$ is symmetric.

1. Consider a distribution $\mu_1$ on the group $X$ of the form     $\mu_1=\gamma_1*m_{\mathbb{Z}(2)}$,
where $\gamma_1$  is a Gaussian distribution on   $\mathbb{R}$ with the characteristic function
$\hat\gamma_1(s)=e^{-\sigma_1 s^2+i\beta_1 s}$, $s\in \mathbb{R}$, where $\sigma_1>0$, $\beta_1\in \mathbb{R}$. Let  $\gamma_j$ be a Gaussian distribution on   $\mathbb{R}$ with the characteristic function
$\hat\gamma_j(s)=e^{-\sigma_j s^2+i\beta_j s}$, $s\in \mathbb{R}$, where $\sigma_j>0$, $\beta_j\in \mathbb{R}$, $j =  2,3, \dots, n$. Assume that conditions    (\ref{11}) hold.
Let $\gamma_j=\lambda_j +\omega_j$,  where $\lambda_j$, $\omega_j$ are measures on $\mathbb{R}$.
Define distributions $\mu_j\in {\rm M}^1(X)$, $j =  2,\dots, n$,  as follows
$$
\mu_j(E\times \{k\})=\begin{cases}\lambda_j(E),   \ \  k=0,\\ \omega_j(E),   \ \  k=1,
\end{cases}
$$
where $E$ is a Borel subset of $\mathbb{R}$. It is obvious that all $\mu_j\in \Lambda$, and    $(\mu_j)_{_{\mathbb{R}}}=\gamma_j$. Taking into account that  $\hat\mu_j(s, 0)=\hat\gamma_j(s)$, $s\in \mathbb{R}$, and  (\ref{09_1}),
  it is easy to see that the characteristic functions  $\hat\mu_j(s, l)$
satisfy Eq. (\ref{05_08_1}).   By Lemma \ref{lem1},  this implies that the conditional distribution of the linear form
$L_2 = b_1\xi_1 + \dots + b_n\xi_n$ given $L_1 = \xi_1 + \dots + \xi_n$ is symmetric.

2. Let $\mu_j\in\Theta$, $j =  1,2, \dots, n$. Assume that the characteristic functions  $\hat\mu_j(s, l)$ are of the form
$$
\hat\mu_j(s, l) = \begin{cases}\exp\{-\sigma_j s^2+i\beta_j s\},
&s\in \mathbb{R}, \ l=0,\\ \kappa_j\exp\{-\sigma'_j s^2+i\beta_j' s\},
&s\in \mathbb{R}, \   l=1.
\\
\end{cases}
$$
Assume that conditions    (\ref{11}) hold and   the equalities
\begin{equation}\label{12_10_17_2}
\sum_{j=1}^nb_j\sigma'_j=0, \quad \sum_{j=1}^nb_j\beta'_j=0
\end{equation}
are fulfilled.  It follows from (\ref{11}) and (\ref{12_10_17_2}) that the characteristic functions  $\hat\mu_j(s, l)$ satisfy Eq. (\ref{05_08_1}). By Lemma \ref{lem1}, this implies    that  the conditional distribution of the linear form
$L_2 = b_1\xi_1 + \dots + b_n\xi_n$ given $L_1 = \xi_1 + \dots + \xi_n$ is symmetric.

Thus, we can not narrow the class of distributions in Theorem \ref{th1} which is characterized by the symmetry of the conditional distribution of one linear form of independent random variables given the other.
\end{remark}
Assume that the number of independent random variables in Theorem  \ref{th1} is equal to 2. It is clear that we may suppose without loss of generality that      $L_1$ and $L_2$   are of the form  $L_1=\xi_1+\xi_2$, $L_2=\xi_1+a\xi_2$, where $a\in {\rm Aut}(X)$, and $a\ne -1$. Discuss  now what happens when
$a=-1$. The following general assertion holds.
\begin{proposition}\label{pr1}   Let $X$ be a locally compact Abelian group, and let   $G=X_{(2)}$. Suppose     $G$  is topologically isomorphic to the group $\mathbb{Z}(2)$. Let $\xi_1$ and $\xi_2$ be independent random variables with values in    $X$ and distributions  $\mu_1$ and $\mu_2$. The conditional distribution of the linear form    $L_2=\xi_1-\xi_2$ given $L_1=\xi_1+\xi_2$ is symmetric if and only if either $\mu_2=\mu_1*\omega_1$ or  $\mu_1=\mu_2*\omega_2$, where $\omega_j\in {\rm M}^1(G)$.
\end{proposition}
{\it Proof}  Proposition \ref{pr1} was proved in \cite{F_solenoid} in the case when the characteristic functions of the distributions  $\mu_j$   do not vanish. The case when the characteristic functions of the distributions  $\mu_j$   can vanish slightly differs from the considered one.

Necessity.   Let $y_0\notin {\overline {Y^{(2)}}}$. Since the character group of the factor-group $Y/{\overline {Y^{(2)}}}$ is topologically isomorphic to the annihilator $A(X, {\overline {Y^{(2)}}})$ and $A(X, {\overline {Y^{(2)}}})=X_{(2)}=G$,
  a decomposition of the group $Y$ with respect to the subgroup ${\overline {Y^{(2)}}}$ is of the form $Y={\overline {Y^{(2)}}}\cup(y_0+{\overline {Y^{(2)}}})$. Since $A(Y, G)={\overline {Y^{(2)}}}$,  it follows from  (\ref{08_1}) that
\begin{equation}\label{08_11_2}
\hat m_{G}(y)=
\begin{cases}
1, & y\in {\overline {Y^{(2)}}},
\\  0, & y\in y_0+{\overline {Y^{(2)}}}.
\end{cases}
\end{equation}
  By Lemma \ref{lem1}, the characteristic functions
$\hat\mu_j(y)$ satisfy equation  $(\ref{42})$ which takes the form
\begin{equation}\label{06_08_17_5}
\hat\mu_1(u+v )\hat\mu_2(u-v)=
\hat\mu_1(u-v)\hat\mu_2(u+v), \quad u, v \in Y.
\end{equation}
Substituting  in $(\ref{06_08_17_5})$ $u=v=y$, we get $\hat\mu_1(2y)=\hat\mu_2(2y)$, $y\in Y$, and hence
\begin{equation}\label{06_08_17_7}
\hat\mu_1(y)=\hat\mu_2(y), \quad y\in {\overline {Y^{(2)}}}.
\end{equation}
Substituting in $(\ref{06_08_17_5})$ $u=y_0+y$, $v=y$, we find
\begin{equation}\label{08_11_1}
\hat\mu_1(y_0+2y)\hat\mu_2(y_0)=
\hat\mu_1(y_0)\hat\mu_2(y_0+2y), \quad y\in Y.
\end{equation}

Assume that there exists an element   $y_0\notin {\overline {Y^{(2)}}}$ such that either $\hat\mu_1(y_0)\ne 0$  or $\hat\mu_2(y_0)\ne 0$.
Suppose for definiteness  that $\hat\mu_1(y_0)\ne 0$. If in this case $\hat\mu_2(y_0)=0$, then (\ref{08_11_1}) implies that $\hat\mu_2(y)=0$ for all $y\in y_0+{\overline {Y^{(2)}}}$. Then it follows from (\ref{08_11_2}) and (\ref{06_08_17_7})  that
$\hat\mu_2(y)=\hat\mu_1(y)\hat m_{G}(y)$, $y\in Y$. Hence    $\mu_2=\mu_1*m_{G}$. Therefore, assume that
 $\hat\mu_2(y_0)\ne 0$.
 Let  ${|\hat\mu_2(y_0)|\le
|\hat\mu_1(y_0)|}$. Put $d_1={\hat\mu_2(y_0))/
\hat\mu_1(y_0)}$. It follows from (\ref{08_11_1})
that $d_1$ is a real number and
\begin{equation}\label{08_08_17_1}
\hat\mu_2(y)=d_1\hat\mu_1(y), \quad y\in y_0+{\overline {Y^{(2)}}}.
\end{equation}
Let $\omega_1$ be a distribution on the group    $G$  with the characteristic function
\begin{equation}\label{08_08_17_2}
\hat\omega_1(y)=
\begin{cases}
1, & \ \  y\in {\overline {Y^{(2)}}},
\\  d_1, & \ \ y\in y_0+{\overline {Y^{(2)}}}.
\end{cases}
\end{equation}
We find from (\ref{06_08_17_7}), (\ref{08_08_17_1}) and (\ref{08_08_17_2})  that
$\hat\mu_2(y)=\hat\mu_1(y)\hat\omega_1(y)$, $y\in Y$.
Hence, $\mu_2=\mu_1*\omega_1$.  If ${|\hat\mu_1(y_0)|\le
|\hat\mu_2(y_0)|}$, we reason similarly, define $d_2$, $\omega_2$ and get that $\mu_1=\mu_2*\omega_2$.

Assume now that there is no an element   $y_0\notin {\overline {Y^{(2)}}}$ such that either $\hat\mu_1(y_0)\ne 0$  or $\hat\mu_2(y_0)\ne 0$. Then   $\hat\mu_1(y)=\hat\mu_2(y)=0$ for all $y\notin {\overline {Y^{(2)}}}$. Taking into account
 (\ref{06_08_17_7}), this implies that $\mu_1=\mu_2$. Moreover, it follows from (\ref{08_11_2}) that $\mu_1=\mu_2*m_{G}$ и $\mu_2=\mu_1*m_{G}$.

Sufficiency implies from Lemma \ref{lem1}. $\blacksquare$

Consider the group $X=\mathbb{R}\times (\mathbb{Z}(2))^2$.  It is interesting to note that      a group analogue of Heyde's theorem on the group $X=\mathbb{R}\times (\mathbb{Z}(2))^2$ may look completely different than on the group    $X=\mathbb{R}\times \mathbb{Z}(2)$. We study the case of two independent random variables and assume that their characteristic functions do not vanish. The following assertion  is true (compare with statement 2 of Theorem \ref{th1}).
 \begin{proposition}\label{pr10.1}  Let $X=\mathbb{R}\times G$, where    $G=(\mathbb{Z}(2))^2$, and let  $\alpha$ be a topological automorphism of the group  $X$ satisfying the condition
\begin{equation}\label{n10_08_1}
{\rm Ker}(I+\alpha)=\{0\}.
\end{equation}
Let $\xi_1$ and $\xi_2$ be
independent random variables with values in   $X$ and distributions
$\mu_1$ and $\mu_2$ with nonvanishing characteristic functions.  If the conditional  distribution of the linear form $L_2 = \xi_1 + \alpha\xi_2$ given $L_1 = \xi_1 + \xi_2$ is symmetric, then $\mu_j\in\Gamma(\mathbb{R})*{\rm M}^1(G)$, $j=1, 2$.
\end{proposition}
Obviously, condition   (\ref{n10_08_1}) implies that $I\pm\alpha\in{\rm Aut}(X)$. For this reason Proposition  \ref{pr10.1} follows from Heyde's theorem for a locally compact Abelian group containing no subgroup topologically isomorphic to the circle group $\mathbb{T}$ (\!\!\cite{Fe20bb}).  Here we give another proof of Proposition \ref{pr10.1}, independent of the results \cite{Fe20bb}. Bearing in mind the description of distributions, which are characterized by the symmetry of the conditional distribution of one linear form of independent random variables with values in a group given the other, the fundamental difference between the groups  $X=\mathbb{R}\times \mathbb{Z}(2)$ and  $X=\mathbb{R}\times (\mathbb{Z}(2))^2$ is that on the group $X=\mathbb{R}\times \mathbb{Z}(2)$ there are  no   topological automorphisms $\alpha$ satisfying   condition  (\ref{n10_08_1}), but on the group $X=\mathbb{R}\times (\mathbb{Z}(2))^2$ such topological automorphisms exist.

{\it Proof of Proposition $\ref{pr10.1}$}      The  group $Y$  is topologically isomorphic to the group $\mathbb{R}\times H$, where $H$ is the character group of the group $G$. The group $H$ is topologically isomorphic to the group $(\mathbb{Z}(2))^2$. Denote by $x=(t, g)$, $t\in \mathbb{R}$, $g\in G$,     elements of the group  $X$,  and by $(s, h)$, $s\in \mathbb{R}$,
 $h\in H$,    elements of the group  $Y$.
Obviously, every automorphism $\alpha\in {\rm Aut}(X)$  is of the form
$\alpha(t, g)=(a t, \alpha_{_G}g)$. In view of (\ref{n10_08_1}) $a\ne -1$ and $\alpha_{_G}g\ne g$ for all $g\in G$, $g\ne 0$. Put $\beta=\tilde\alpha_{_G}$. It is obvious that $\beta h\ne h$ for all $h\in H$, $h\ne 0$.  Denote by   $p$, $q$, $p+q$   nonzero elements of the group $H$. It is clear that without loss of generality we may suppose that
\begin{equation}\label{n11_08_3}
\beta p=q, \ \beta q=p+q.
\end{equation}
By Lemma \ref{lem1}, the characteristic functions  $\hat\mu_j(s, h)$ satisfy Eq.  (\ref{42})  which takes the form
\begin{equation}\label{n10_08_2}
\hat\mu_1(s_1+s_2, h_1+h_2)\hat\mu_2(s_1+as_2, h_1+\beta h_2)$$$$=
\hat\mu_1(s_1-s_2, h_1+h_2)\hat\mu_2(s_1-as_2, h_1+\beta h_2), \quad s_1, s_2\in \mathbb{R}, \ h_1, h_2\in H.
\end{equation}
Substituting  $h_1=h_2=0$ in (\ref{n10_08_2})  we obtain
\begin{equation}\label{n10_08_3}
\hat\mu_1(s_1+s_2, 0)\hat\mu_2(s_1+as_2, 0)$$$$=
\hat\mu_1(s_1-s_2, 0)\hat\mu_2(s_1-as_2, 0), \quad s_1, s_2\in \mathbb{R}.
\end{equation}
Hence, taking into account Lemma \ref{lem1}, it follows from Heyde's theorem that
\begin{equation}\label{n10_08_4}
\hat\mu_j(s, 0)=e^{-\sigma_j s^2+i\beta_j s}, \quad s\in \mathbb{R}, \   j = 1, 2,
\end{equation}
where $\sigma_j\ge 0$,  $\beta_j\in \mathbb{R}$. Substituting (\ref{n10_08_4}) into  (\ref{n10_08_3}) we get
$\beta_1+a\beta_2=0$. For this reason replacing the distributions $\mu_j$ with their shifts $\mu_j*E_{-\beta_j}$, we can suppose that $\beta_1=\beta_2=0$ in (\ref{n10_08_4}), i.e.
\begin{equation}\label{n10_08_5}
\hat\mu_j(s, 0)=e^{-\sigma_j s^2}, \quad s\in \mathbb{R}, \   j = 1, 2,
\end{equation}
is fulfilled from the very beginning. By Lemma \ref{lem7}, for each fixed $h\in H$ the functions $\hat\mu_j(s, h)$ can be extended to the complex plane as entire functions in   $s$. It is easy to see that Eq.  (\ref{n10_08_2})  holds for all complex $s_j$ and $h_j\in H$. It also follows from (\ref{n10_08_2})  that the functions  $\hat\mu_j(s, h)$ do not vanish in the complex plane.  Hence, we have the representation
\begin{equation}\label{n11_08_1}
\hat\mu_j(s, h)=\exp\{\psi_{j}(s,h)\}, \quad s\in \mathbb{R}, \quad h\in H, \quad j=1, 2,
\end{equation}
where $\psi_{j}(s,h)$ are entire functions in $s$.
By Lemma \ref{lem7}, (\ref{n10_08_5}) implies that for each fixed   $h\in H$ the functions $\hat\mu_j(s, h)$ are entire functions in    $s$  of the order at most   2.  For this reason it follows from (\ref{n11_08_1}) that
\begin{equation}\label{n10_08_6}
\psi_{j}(s,h)=a_{jh}s^2+b_{jh}s+c_{jh}, \quad s\in \mathbb{R}, \quad h\in H, \quad j=1, 2,
\end{equation}
where $a_{jh}$, $b_{jh}$, $c_{jh}$ are some complex constants.  Substituting (\ref{n11_08_1}) into (\ref{n10_08_2})  and taking into account (\ref{n10_08_6})  from the resulting equation, we find
\begin{equation}\label{n11_08_2}
a_{1,h_1+h_2}(s_1+s_2)^2+a_{2,h_1+\beta h_2}(s_1+as_2)^2$$$$=a_{1,h_1+h_2}(s_1-s_2)^2+a_{2,h_1+\beta h_2}(s_1+as_2)^2, \quad s_1, s_2\in \mathbb{R}, \quad h_1, h_2\in H.
\end{equation}
Substituting  in succession  $h_1=h_2=0$, then $h_1=h_2=p$, then  $h_1=h_2=q$, then  $h_1=h_2=p+q$ in (\ref{n11_08_2}). Taking into account (\ref{n11_08_3}) and equating in the left and right sides of each of the   obtained equations the coefficients at $s_1s_2$, we get
\begin{equation}\label{n11_08_4}
a_{1,0}+aa_{2,0}=0, \ a_{1,0}+aa_{2,p+q}=0, \ a_{1,0}+aa_{2,p}=0, \ a_{1,0}+aa_{2,q}=0.
\end{equation}
In view of (\ref{n10_08_5}), it follows from  (\ref{n11_08_4}) that
\begin{equation}\label{n11_08_5}
a_{2,h}=-\sigma_2, \quad h\in H.
\end{equation}
Arguing similarly,   substitute  in succession     $h_1=h_2=0$, then $h_1=p$, $h_2=p+q$, then $h_1=q$, $h_2=p$, then  $h_1=p+q$, $h_2=q$  in (\ref{n11_08_2}). Taking into account (\ref{n10_08_5}), we find
\begin{equation}\label{n11_08_6}
a_{1,h}=-\sigma_1, \quad h\in H.
\end{equation}
Substituting (\ref{n11_08_1}) into (\ref{n10_08_2})  and taking into account (\ref{n10_08_6}), from the resulting equation, we find
\begin{equation}\label{n11_08_7}
b_{1,h_1+h_2}(s_1+s_2)+b_{2,h_1+\beta h_2}(s_1+as_2)$$$$=b_{1,h_1+h_2}(s_1-s_2)+b_{2,h_1+\beta h_2}(s_1-as_2), \quad s_1, s_2\in \mathbb{R}, \quad h_1, h_2\in H.
\end{equation}
Reasoning as above, equating in the left and right sides of    Eq. (\ref{n11_08_7}) the coefficients at $s_2$, and taking into account (\ref{n10_08_5}), we find
\begin{equation}\label{n11_08_8}
b_{1,h}=b_{2,h}=0, \quad h\in H.
\end{equation}
From (\ref{n11_08_5}), (\ref{n11_08_6}) and (\ref{n11_08_8}) we get that $\psi_{j}(s,h)=-\sigma_js^2+c_{jh}$, $j=1, 2$. Hence,
\begin{equation}\label{n11_08_9}
\hat\mu_j(s, h)=\exp\{-\sigma_js^2+c_{jh}\}, \quad s\in \mathbb{R}, \quad h\in H, \quad j=1, 2.
\end{equation}

 Denote by $\gamma_j$  the Gaussian distribution on   $\mathbb{R}$ with the characteristic function $\hat\gamma_j(s)=\exp\{-\sigma_js^2\}$, and by $\omega_j$ the distribution on the group $G$ with characteristic function
\begin{equation}\label{n11_08_10}
\hat\omega_j(h)=\exp\{c_{jh}\}, \quad  h\in H, \quad j=1, 2.
\end{equation}
It follows from  (\ref{n11_08_9}) and (\ref{n11_08_10}) that
$\hat\mu_j(s, h)=\hat\gamma_j(s)\hat\omega_j(h)$, $s\in \mathbb{R}$, $h\in H$. Hence,  $\mu_j=\gamma_j*\omega_j$,  $j=1, 2$, which proves the proposition.  $\blacksquare$

\begin{remark}\label{re2}     Heyde's theorem is similar to the well-known Skitovich--Darmois theorem, where the class of Gaussian distributions on the real line is characterized by the independence of two linear forms of independent random variables  (see e.g. \cite[Chapter 3]{KaLiRa}).
\begin{a*} Let $\xi_j$, $j = 1, 2,\dots, n$, $n \ge 2,$ be independent random variables with distributions  $\mu_j$. Let $a_j$, $b_j$ be nonzero real numbers.
Assume that the linear forms $L_1=a_1\xi_1+\cdots+a_n\xi_n$ and
$L_2=b_1\xi_1+\cdots+b_n\xi_n$ are independent. Then  all distributions $\mu_{j}$   are Gaussian, possibly degenerate.
\end{a*}
Compare Theorem \ref{th1} with an   analogue of the Skitovich--Darmois theorem for the group   $X=\mathbb{R}\times \mathbb{Z}(2)$. It is not difficult to prove that the following assertion holds.
\begin{b*} Let $X=\mathbb{R}\times \mathbb{Z}(2)$, and let $a_j$, $b_j$, $j = 1, 2,\dots, n$, $n \ge 2,$  be topological automorphisms of  the group $X$. Let $\xi_j$ be independent random variables with values in     $X$ and distributions  $\mu_j$.
Assume that the linear forms $L_1=a_1\xi_1+\cdots+a_n\xi_n$ and
$L_2=b_1\xi_1+\cdots+b_n\xi_n$ are independent. Then  all distributions $\mu_{j}$   are Gaussian.
\end{b*}
We see that  on the one hand,    the class of Gaussian distributions on the real line is characterized  both by the independence of two linear forms of independent random variables and by the  symmetry of the conditional distribution of one linear form given the other. On the other hand,    on the group $X=\mathbb{R}\times \mathbb{Z}(2)$ the independence of two linear forms of independent random variables  characterizes only Gaussian distributions, whereas the  symmetry of the conditional distribution of one linear form given the other   characterizes much wider class of distributions    (see Theorem \ref{th1}).
\end{remark}

\section{    Heyde's characterization theorem for the group $\mathbb{R}\times\mathbb{T}$}

Let $X=\mathbb{R}\times\mathbb{T}$. Denote by $x=(t, z)$, $t\in\mathbb{R}$, $z\in\mathbb{T}$,  elements of the group $X$.    The   group $Y$  is topologically isomorphic to the group $\mathbb{R}\times\mathbb{Z}$.
Denote by $y=(s,n)$, $s\in\mathbb{R}$, $n\in\mathbb{Z}$, elements of the group $Y$.
 It is easy to verify that each   automorphism $\alpha\in{\rm
Aut}(X)$ is defined
by a matrix
$\left(\begin{matrix}a&b\\ 0&\pm 1\end{matrix}\right)$, where
$a, b\in \mathbb{R}$, $a\neq 0$, and $\alpha$ is given by the formula
$$
\alpha(t,z)=(at,e^{ibt}z^{\pm 1}),\quad
t\in\mathbb{R}, \ z\in\mathbb{T}.
$$
Then the adjoint automorphism
$\tilde\alpha\in{\rm
Aut}(Y)$ is of the form
$$
\tilde\alpha(s,n)=(as+bn,\pm n),\quad
s\in\mathbb{R}, \ n\in\mathbb{Z}.
$$
We will identify the automorphisms ${\alpha}$ and ${\tilde\alpha}$ with
with the corresponding
matrix $\left(\begin{matrix}a&b\\ 0&\pm 1\end{matrix}\right).$

Denote by $G$ the subgroup of $X$ generated by the element of order 2.  Let $\alpha\in{\rm
Aut}(X)$. Assume that   $\alpha$ is of the form $\alpha=\left(\begin{matrix}a&b\\ 0& 1\end{matrix}\right)$, where   $a\ne 1$. Define the continuous monomorphism   $\tau\colon\mathbb{R}\times\mathbb{Z}(2)\to X$ by the formula \begin{equation}\label{22_1}
\tau(t, k)=\left(t, (-1)^k e^{ibt\over {a-1}}\right), \quad t\in\mathbb{R}, \ k\in \mathbb{Z}(2).
\end{equation}
Put $F=\tau(\mathbb{R})$. Then $F$ is a close subgroup of   $X$ and $F$ is topologically isomorphic to   $\mathbb{R}$. Note that $\tau(\mathbb{Z}(2))=G$.
We retain the notation $G$, $\tau$ and $F$  throughout this section.

In this    section we study an  analogue of Heyde's theorem
for  two  independent random variables with values
in the group $X$. Without loss of generality, we can assume that   the linear forms  $L_1$ and $L_2$   are    $L_1=\xi_1+\xi_2$, $L_2=\xi_1+\alpha\xi_2$, where $\alpha\in {\rm Aut}(X)$.
The following statement can be viewed as an  analogue of Heyde's theorem
for two independent random variables with values
in the group $X$ in the case when $\alpha=\left(\begin{matrix}a&b\\ 0& 1\end{matrix}\right)$. It turns out that in this case the description of possible distributions that are characterized by the symmetry of the conditional distribution of the linear form
$L_2=\xi_1+\alpha\xi_2$
 given $L_1 = \xi_1 + \xi_2$ depends only on $a$.
\begin{theorem}\label{th2}
Let $X=\mathbb{R}\times\mathbb{T}$,  and  let  $\alpha=\left(\begin{matrix}a&b\\ 0& 1\end{matrix}\right)\in{\rm
Aut}(X)$.
   Let $\xi_1$ and $\xi_2$ be independent random variables with values in the group
       $X$  and distributions
  $\mu_1$ and $\mu_2$.
If the conditional distribution of the linear form
$L_2=\xi_1+\alpha\xi_2$
 given $L_1 = \xi_1 + \xi_2$ is symmetric, then the following statements hold.

$1$. If $a<0$ and $a\ne - 1$, then  there is the following alternative.

$1a$. The  distributions $\mu_{j}$ are represented in the form $\mu_j=\tau(M_j)*E_{x_j}$,
  where $x_j\in X$, and  $M_j$ are distributions on the group  $\mathbb{R}\times\mathbb{Z}(2)$   belonging  to the class $\Lambda$, and at least one of the distributions  $\mu_j$ is represented in the form $\mu_j=\gamma_j*m_G*E_{x_j},$ where  $\gamma_j\in \Gamma(F)$.

$1b$. The distributions $\mu_{j}$ are represented in the form $\mu_j=\tau(M_j)*E_{x_j}$,
  where $x_j\in X$, and $M_j$ are distributions on the group  $\mathbb{R}\times\mathbb{Z}(2)$   belonging to the class $\Theta$, and  the characteristic functions  $\hat\mu_j(y)$ do not vanish.

$2$. If $a=-1$, then   the distributions $\mu_j$  are of the form $\mu_j=\lambda_j*E_{x_j},$
where $x_j\in X$, and $\lambda_j\in{\rm M}^1(F\times G)$. Furthermore	
either $\lambda_2=\lambda_1*\omega_1$    or $\lambda_1=\lambda_2*\omega_2$, where $\omega_j\in {\rm M}^1(G)$.

$3$.  If   $a> 0$, then some shifts of the distributions $\mu_j$
are supported in $G$.
\end{theorem}

The proof of Theorem \ref{th2} is based on the following statement.
\begin{lemma}\label{lem23.1}   Let $X=\mathbb{R}\times\mathbb{T}$,    and let   $\alpha=\left(\begin{matrix}a&b\\ 0& 1\end{matrix}\right)\in{\rm
Aut}(X)$.
   Let $\xi_1$ and $\xi_2$ be independent random variables with values in the group
       $X$  and distributions
  $\mu_1$ and $\mu_2$.
Assume that the conditional distribution of the linear form
$L_2=\xi_1+\alpha\xi_2$
 given $L_1 = \xi_1 + \xi_2$ is symmetric.  Then there exist     elements $x_j\in X$ such that if $a\ne 1$, then the  distributions $\lambda_j=\mu_j*E_{-x_j}$   are supported in  the subgroup  $F\times G$, and if $a=1$, then $\lambda_j$ are supported in the subgroup  $G$. Moreover, if $\eta_j$ are independent random variables with values in the group
       $X$  and distributions
  $\lambda_j$, then   the conditional distribution of the linear form   $N_2=\eta_1+\alpha\eta_2$ given $N_1=\eta_1+\eta_2$  is symmetric.
\end{lemma}
{\it Proof} \ Put
 $L={\rm Ker}(I-\tilde\alpha)$. It is obvious that $L=\{(s, n):(1-a)s-bn=0\}$. Put $K=A\left(X, L^{(2)}\right)$.
 By Lemma \ref{lem1},
 the characteristic functions
 $\hat\mu_j(y)$ satisfy Eq.
(\ref{42}) which takes the form
\begin{equation}\label{05_11_1}
\hat\mu_1(u+v )\hat\mu_2(u+\tilde\alpha v )=
\hat\mu_1(u-v )\hat\mu_2(u-\tilde\alpha v), \quad u, v \in Y.
\end{equation}
 Consider the restriction of Eq.  (\ref{05_11_1}) to the subgroup $L$. Since $\tilde\alpha y=y$ for all $y\in L$, we have
\begin{equation}\label{05_11_2}
\hat\mu_1(u+v )\hat\mu_2(u+ v )=
\hat\mu_1(u-v )\hat\mu_2(u- v), \quad u, v \in L.
\end{equation}
Substituting $u=v=y$ in (\ref{05_11_2}), we obtain
\begin{equation}\label{05_11_3}
\hat\mu_1(y)\hat\mu_2(y)=1, \quad y \in L^{(2)}.
\end{equation}
Hence, $|\hat\mu_1(y)|=|\hat\mu_2(y)|=1$, $y \in L^{(2)}$. It follows from this  that there exist elements $p_j\in X$ such that
\begin{equation}\label{27_01_1}
\hat\mu_j(y)=(p_j, y), \ y\in L^{(2)}, \ j=1, 2.
\end{equation}
Substituting (\ref{27_01_1}) in (\ref{05_11_3}), we receive that $p_1+p_2\in K$.
Let $p_j=(t_j, z_j)$, $j=1, 2$.

1. Suppose that $a\ne 1$.   We have
$$L^{(2)}=\left\{\left({2bn\over {1-a}}, 2n\right): n\in \mathbb{Z}\right\},$$
and hence,
$$
K=A\left(X, L^{(2)}\right)=\left\{\left(t, \pm e^{ibt\over {a-1}}\right): t\in \mathbb{R}\right\}.
$$
 It is clear that
$K=F\times G$ and $\tau$ is a topological isomorphism of the groups    $\mathbb{R}\times \mathbb{Z}(2)$
and $K$.

Put $\tilde x=\left(t_1,  e^{{ibt_1\over {a-1}}}\right)\in K$. Verify that  $x_1=p_1-\tilde x$ and $x_2=-p_1+\tilde x$ are the required elements. Consider the distributions $\lambda_j=\mu_j*E_{-x_j}$.   It follows from $p_1+p_2, \tilde x\in K$ that $\hat\lambda_1(y)=\hat\lambda_2(y)=1$ for all $y\in L^{(2)}$. By Lemma \ref{lem2}, this implies that
$\sigma(\lambda_j)\subset A\left(X, L^{(2)}\right)=K$, $j=1, 2$. It is obvious that $x_j\in \mathbb{T}$. Note that the restriction of the automorphism $\alpha$ to the subgroup $\mathbb{T}$ is the identity automorphism.  Hence,
$x_1+\alpha x_2=0$. This implies that the functions $(-x_1, y)$ and $(-x_2, y)$  satisfy Eq.   (\ref{05_11_1}). Since    the characteristic functions $\hat\mu_j(y)$  satisfy Eq.   (\ref{05_11_1}), we obtain that the characteristic functions
  $\hat\lambda_j(y)$ also satisfy Eq.  (\ref{05_11_1}). Let
$\eta_j$ be independent random variables with values in the group
 $X$ and distributions  $\lambda_j$.
Then by Lemma \ref{lem1}, the conditional distribution of the linear form $N_2 = \eta_1 +\alpha\eta_2$ given
 $N_1 = \eta_1 + \eta_2$ is symmetric.

2. Let $a=1$, $b=0$, i.e.  $\alpha=I$. Hence, $L=Y$, and then $K=A\left(X, L^{(2)}\right)=G$. It is easy to see that we can put $x_j=p_j$, $j=1, 2$.

3. Let $a=1$, $b\ne 0$. Then  $L=\{(s, 0): s\in \mathbb{R}\}=\mathbb{R}$, and hence, $K=A(X, \mathbb{R})=\mathbb{T}$. This implies that $t_1+t_2=0$.
Put $\hat x=(0, z_1z_2e^{-ibt_1})$ and  $x_1=p_1-\hat x$, $x_2=p_2$. Consider the distributions $\lambda_j=\mu_j*E_{-x_j}$.
Since $\hat x\in \mathbb{T}$, we obtain from (\ref{27_01_1}) that $\hat\lambda_1(s, 0)=1$ for all $s\in \mathbb{R}$. By Lemma \ref{lem2}, it follows from this that
$\sigma(\lambda_j)\subset A\left(X, \mathbb{R}\right)=\mathbb{T}$, $j=1, 2$.
  It is obvious that $x_1+\alpha x_2=0$. Taking this into account, and   taking into account the fact that the characteristic functions  $\hat\mu_j(y)$  satisfy Eq. (\ref{05_11_1}), we get that
the characteristic functions   $\hat\lambda_j(y)$  also satisfy Eq.  (\ref{05_11_1}).  Let
$\eta_j$ be independent random variables with values in the group
 $X$ and distributions  $\lambda_j$. By Lemma \ref{lem1},  the conditional distribution of the linear form $N_2 = \eta_1 + \alpha\eta_2$ given
  $N_1 = \eta_1 + \eta_2$ is symmetric.  Taking into account that the restriction of the automorphism  $\alpha$ to the subgroup $\mathbb{T}$
 is the identity automorphism, we reduce  the proof to the case when
 $X=\mathbb{T}$ and $\alpha=I$. To complete the proof we argue in the same way as in case 2.  $\blacksquare$
\medskip

{\it Proof of Theorem $\ref{th2}$}  \  Assume that $a\ne 1$.  We apply Lemma \ref{lem23.1} and get independent random variables $\eta_j$ with values in the group $X$ and distributions $\lambda_j$. Put $K=F\times G$. Since $\sigma(\lambda_j)\subset K$,  we may consider $\eta_j$  as independent random variables with values in the group $K$. It is easy to see that the restriction of the   automorphism $\alpha$ to the subgroup $K$ is a topological automorphism of the subgroup $K$, and $\alpha$ operates
 on the elements of  $K$ as follows
\begin{equation}\label{11_0}
\alpha\left(t, \pm e^{ibt\over {a-1}}\right)=\left(at, \pm e^{ibat\over {a-1}}\right), \quad t\in \mathbb{R}.
\end{equation}
 We also note that since $\tau$ is a topological isomorphism of the groups  $\mathbb{R}\times \mathbb{Z}(2)$
and $K$. Moreover, $\tau$ generates in a natural way an isomorphism of the convolution semigroups of probability distributions ${{\rm M}^1(\mathbb{R}\times \mathbb{Z}(2))}$ and ${\rm M}^1(K)$. We keep the notation $\tau$ for this isomorphism of the convolution semigroups   ${{\rm M}^1(\mathbb{R}\times \mathbb{Z}(2))}$ and ${\rm M}^1(K)$.

Let us prove statement 1. Put $\zeta_j=\tau^{-1} \eta_j$. Then $\zeta_j$ are independent random variables   with values in the group
$\mathbb{R}\times \mathbb{Z}(2)$. Denote by $M_j$ the distribution of the random variable  $\zeta_j$. Then
$\lambda_j=\tau(M_j)$, $j=1, 2$. It follows from (\ref{11_0}) that  $\tau^{-1}\alpha=a\tau^{-1}$, and hence
$\tau^{-1}\alpha\eta_j=a\zeta_j.$
The symmetry of the conditional distribution of the linear form $N_2=\eta_1+\alpha\eta_2$ given $N_1=\eta_1+\eta_2$ implies the symmetry of the conditional distribution of the linear form $P_2 = \zeta_1 + a\zeta_2$  given $P_1 = \zeta_1 +
\zeta_2$. Since $\tau(\Gamma(\mathbb{R}))=\Gamma(F)$ and $\tau(m_{\mathbb{Z}(2)})=m_G$, statement 1 follows from Theorem \ref{th1}.

Let us prove  statement $2$.   Since $a=-1$, it follows from (\ref{11_0}) that  the restriction of the automorphism $\alpha$ to the subgroup $K$ coincides with $-I$. It means that if we consider $\eta_j$ as independent random variables with values in the group   $K$, then the conditional distribution of the linear form $N_2 = \eta_1 -\eta_2$ given $N_1 = \eta_1 +\eta_2$ is symmetric. Applying Proposition \ref{pr1} to the group $K$, to  the independent random variables $\eta_1$ and $\eta_2$, and to the linear forms $N_1$ and $N_2$,   we prove statement 2.

Let us prove  statement 3.  It follows from Lemma \ref{lem23.1} that it suffices to consider the case when   $a\ne 1$. We retain the notation used in the proof of statement 1.   Applying Corollary  \ref{co1} to the random variables $\zeta_j$ and linear forms   $P_1$  and $P_2$, we obtain that some shifts of the distributions  $M_j$ are supported in the subgroup
 $\mathbb{Z}(2)$. Since the convolution semigroups  ${\rm M}^1(\mathbb{R}\times \mathbb{Z}(2))$ and ${\rm M}^1(K)$
are isomorphic and $\tau(\mathbb{Z}(2))=G$, some shifts of the distributions   $\lambda_j$
are supported in the subgroup $G$. Hence, some shifts of the distributions    $\mu_j$ are also supported in   $G$. $\blacksquare$
\begin{remark}\label{re3}  Let $X=\mathbb{R}\times\mathbb{T}$. Suppose that $\alpha=\left(\begin{matrix}a&b\\ 0& 1\end{matrix}\right)\in {\rm Aut}(X)$, where    $a$ is such as in each of statements 1--3 of Theorem \ref{th2}. Then we can construct distributions $\mu_1$ and $\mu_2$ such as in  the corresponding statement of  Theorem \ref{th2} with the following property:
 if
 $\xi_1$ and  $\xi_2$    are independent random variables with values in the group
       $X$  and distributions
 $\mu_1$ and $\mu_2$,  then the conditional distribution of the linear form
$L_2 = \xi_1 + \alpha\xi_2$
 given $L_1 = \xi_1 + \xi_2$ is symmetric.
Thus, we can not narrow the class of distributions in Theorem \ref{th2} which is characterized by the symmetry of the conditional distribution of one linear form of independent random variables given the other.
  \end{remark}
\begin{remark}\label{re4}  Let $X=\mathbb{R}\times\mathbb{T}$, and let the continuous homomorphism ${\tau\colon\mathbb{R}\times\mathbb{Z}(2)\to X}$ be defined by formula (\ref{22_1}).  Then the adjoint homomorphism    ${\tilde\tau\colon\mathbb{R}\times\mathbb{Z}\to \mathbb{R}\times\mathbb{Z}(2)}$ is of the form $\tilde\tau(s, n)=\left(s+{bn\over{a-1}}, n\ ({\rm mod} \ 2)\right)$.
Assume that in statement $1b$ of Theorem \ref{th2} the distributions $\mu_j$ are represented in the form $\mu_j=\tau(M_j)$,   where
$M_j\in\Theta$,  and the characteristic functions    $\hat M_j(s, l)$ are of the form
$$
\hat M_j(s, l) = \begin{cases}e^{-\sigma_j s^2+i\beta_j s}, \ \  &s\in \mathbb{R}, \ l=0,\\ \kappa_je^{-\sigma_j' s^2+i\beta_j's},
&s\in \mathbb{R}, \   l=1,
\\
\end{cases}
$$
then the characteristic functions
$\hat\mu_j(s, n)$ can be written as follows
$$
\hat\mu_j(s, n) = \begin{cases}\exp\left\{-\sigma_j \left(s+{bn\over{a-1}}\right)^2+i\beta_j \left(s+{bn\over{a-1}}\right)\right\}, \ \ \  &s\in \mathbb{R}, \ n\in \mathbb{Z}^{(2)},\\ \kappa_j\exp\left\{-\sigma_j' \left(s+{bn\over{a-1}}\right)^2+i\beta_j'\left(s+{bn\over{a-1}}\right)\right\},
 &s\in \mathbb{R}, \   n\in \mathbb{Z}\backslash\mathbb{Z}^{(2)} \ .
\\
\end{cases}
$$
 \end{remark}
Let $X=\mathbb{R}\times\mathbb{T}$, and let $\alpha=\left(\begin{matrix}a&b\\ 0& 1\end{matrix}\right)\in{\rm
Aut}(X)$. Theorem \ref{th2} gives us a complete description of possible distributions of independent random variables $\xi_1$ and $\xi_2$ with values in the group $X$ on the assumption that the conditional distribution of the linear form
$L_2=\xi_1+\alpha\xi_2$
 given $L_1 = \xi_1 + \xi_2$ is symmetric.

 Now we will discuss what happens when $\alpha=\left(\begin{matrix}a&b\\ 0& -1\end{matrix}\right)\in{\rm
Aut}(X)$.

 Let either  $a>0$  or $a=-1$ and $b=0$. Assume that $\xi_1$ and  $\xi_2$    are independent random variables with values in the group
       $X$  and distributions
 $\mu_1$ and $\mu_2$. Suppose  that the conditional distribution of the linear form
$L_2 = \xi_1 + \alpha\xi_2$
 given $L_1 = \xi_1 + \xi_2$ is symmetric. What   can we say about the distributions $\mu_j$? By Lemma \ref{lem1},
 the characteristic functions $\hat\mu_j(y)$ satisfy Eq.
(\ref{05_11_1}), which takes the form
\begin{equation}\label{08_11_3}
\hat\mu_1(s_1+s_2, n_1+n_2)\hat\mu_2(s_1+as_2+bn_2,n_1-n_2)$$$$=
\hat\mu_1(s_1-s_2, n_1-n_2)\hat\mu_2(s_1-as_2-bn_2,n_1+n_2), \quad s_1, s_2\in\mathbb{R}, \ n_1, n_2\in\mathbb{Z}.
\end{equation}

1. Let $a>0$.
Substitute  $n_1=n_2=0$ in (\ref{08_11_3}). It follows from Lemma \ref{lem1} and   Heyde's theorem  that
$$
\hat\mu_j(s, 0)=\exp\{-\sigma_j s^2+i\beta_j s\},\quad s\in \mathbb{R}, \ j=1, 2,
$$
where  $\sigma_j\ge 0$, $\beta_j\in \mathbb{R}$. Furthermore,  $\sigma_1+a\sigma_2=0$ and $\beta_1+a\beta_2=0$. Since $a>0$, we have $\sigma_1=\sigma_2=0$. Thus,
$$
\hat\mu_j(s, 0)=\exp\{i\beta_j s\},\quad s\in\mathbb{R}, \ j=1, 2.
$$
Put
$\lambda_j=\mu_j*E_{(-\beta_j, 1)}$. Since the characteristic functions  $\hat\mu_j(s, n)$  satisfy Eq. (\ref{08_11_3}) and $\beta_1+a\beta_2=0$, the characteristic functions $\hat\lambda_j(s, n)$ also satisfy Eq.   (\ref{08_11_3}). Let $\eta_j$  be independent random variables with values in the group  $X$ and distributions $\lambda_j$. By Lemma
 \ref{lem1}, the conditional distribution of the linear form  $N_2=\eta_1+\alpha\eta_2$ given $N_1=\eta_1+\eta_2$ is symmetric.
 It is obvious that $\hat\lambda_1(s, 0)=\hat\lambda_2(s, 0)=1$ for all $s\in \mathbb{R}$. By Lemma \ref{lem2}, it follows from this that
$\sigma(\lambda_j)\subset A(X, \mathbb{R})=\mathbb{T}$, $j=1, 2$. Note that the restriction of the automorphism $\alpha$ to the subgroup $\mathbb{T}$ coincides with $-I$. Hence, if we consider $\eta_j$ as independent random variables with values in the circle group   $\mathbb{T}$ and distributions $\lambda_j$, then the conditional distribution of the linear form $N_2 = \eta_1 -\eta_2$ given $N_1 = \eta_1 +\eta_2$ is symmetric.
Applying Proposition \ref{pr1} to the circle group $\mathbb{T}$, to the independent random variables $\eta_1$ and  $\eta_2$, and to the linear forms $N_1$ and $N_2$,  we obtain that either $\lambda_2=\lambda_1*\omega_1$    or $\lambda_1=\lambda_2*\omega_2$, where $\omega_j\in{\rm M}^1(G)$.

2.  Let $a=-1$ and $b=0$. Then $\alpha=-I$, and the description of the distributions $\mu_j$ follows from Proposition \ref{pr1}.

We see that in cases  1 and   2 the description of possible distributions of independent random variables $\xi_j$ which are characterized by the symmetry of the conditional distribution of the linear form   $L_2=\xi_1+\alpha\xi_2$ given $L_1=\xi_1+\xi_2$ is reduced to the case when  $L_2=\xi_1-\xi_2$, i.e. to Proposition  \ref{pr1}.

Now suppose   that either $a=-1$ and $b\ne 0$  or $a<0$ and $a\ne -1$. In each of these cases we construct distributions  $\mu_1$ and $\mu_2$ with the following property:    if $\xi_1$ and
$\xi_2$ are independent random variables with values in the group
       $X$  and distributions
 $\mu_1$ and $\mu_2$, then the conditional distribution of the linear form
$L_2 = \xi_1 + \alpha\xi_2$
 given $L_1 = \xi_1 + \xi_2$ is symmetric.

3. Let $a=-1$ and $b\ne 0$. In this case we construct   two different types of distributions  $\mu_j$.

$3a$. Consider the subgroup  $H={\rm Ker}(I+\alpha)$. Then $H=\left\{\left({2\pi n\over b}, z\right): n\in \mathbb{Z}, z\in \mathbb{T}\right\}$ and $H$ is topologically isomorphic to the group $\mathbb{Z}\times\mathbb{T}$.
Let either $\mu_2=\mu_1*\omega_1$, where $\mu_1\in {\rm M}^1(H)$, $\omega_1\in {\rm M}^1(G)$ or  $\mu_1=\mu_2*\omega_2$, where   $\mu_2\in {\rm M}^1(H)$, $\omega_2\in {\rm M}^1(G)$.
By Proposition \ref{pr1}, if $\xi_1$ and
$\xi_2$ are independent random variables with values in the group
       $H$  and distributions
 $\mu_1$ and $\mu_2$, then the conditional distribution of the linear form
$L_2 = \xi_1
-\xi_2$
 given $L_1 = \xi_1 + \xi_2$ is symmetric. Note that the restriction of  the automorphism  $\alpha$ to the subgroup $H$ coincides with $-I$. Hence, if we consider $\xi_1$ and
$\xi_2$   as   independent random variables with values in the group   $X$, then the conditional distribution of the linear form $L_2 = \xi_1 +
\alpha\xi_2$  given $L_1 = \xi_1 + \xi_2$ is also symmetric.

$3b$. Let  $\{c_n\}_{n\in \mathbb{Z}}$ be a sequence of complex numbers satisfying the conditions
$$
c_0=1, \ \ c_{-n}=\bar c_n,  \  n\in \mathbb{Z},  \ \ \sum_{n\in \mathbb{Z}, \ n\ne 0}|c_n|<1.
$$
It follows from this that
$$
\rho(z)=\sum_{n\in \mathbb{Z}}c_nz^n>0, \ z\in \mathbb{T}, \quad \int_\mathbb{T}\rho(z)dm_\mathbb{T}(z)=1.
$$
Assume that   $c_{2n}=0$ for all $n\in \mathbb{Z}$, $n\ne 0$. Denote by $\lambda$ the distribution on the circle group $\mathbb{T}$ with the density $
\rho(z)$ with respect to $m_\mathbb{T}$. Then $\hat\lambda(2n)=0$ for all $n\in \mathbb{Z}$, $n\ne 0$.
Let $\mu$ be an arbitrary distribution on the group $\mathbb{R}$. Consider the groups $\mathbb{R}$ and   $\mathbb{T}$ as subgroups of $X$ and put $\mu_1=\mu*m_\mathbb{T}$, $\mu_2=\mu*\lambda$.
We have then
\begin{equation}\label{08_01_18_1}
\hat\mu_1(s, 0)=\hat\mu_2(s, 0)=\hat\mu(s), \ s\in\mathbb{R},
\end{equation}
\begin{equation}\label{08_01_18_2}
\hat\mu_1(s, n)=\hat\mu_2(s, 2n)=0, \ s\in\mathbb{R}, \ n\in \mathbb{Z}, \ n\ne 0.
\end{equation}
It is not difficult to check that the characteristic functions $\hat\mu_j(s, n)$ satisfy the equation
\begin{equation}\label{08_01_18_5}
\hat\mu_1(s_1+s_2, n_1+n_2)\hat\mu_2(s_1-s_2+bn_2,n_1-n_2)$$$$=
\hat\mu_1(s_1-s_2, n_1-n_2)\hat\mu_2(s_1+s_2-bn_2,n_1+n_2), \quad  s_1, s_2\in\mathbb{R}, \ n_1, n_2\in\mathbb{Z}.
\end{equation}
Indeed, if $n_1=n_2=0$, then taking into account (\ref{08_01_18_1}), (\ref{08_01_18_5})  turns into the equality
$$
\hat\mu(s_1+s_2)\hat\mu(s_1-s_2)=
\hat\mu(s_1-s_2)\hat\mu_(s_1+s_2), \quad  s_1, s_2\in\mathbb{R}.
$$
If either $n_1\ne 0$ or $n_2\ne 0$, then   (\ref{08_01_18_2})  implies that both parts of Eq.  (\ref{08_01_18_5}) are equal to zero.
Let $\xi_1$ and
$\xi_2$ be independent random variables with values in the group
       $X$ and distributions
 $\mu_1$ and $\mu_2$.
   By Lemma \ref{lem1}, it follows from (\ref{08_01_18_5}) that the conditional distribution of the linear form $L_2 = \xi_1 +
\alpha\xi_2$  given $L_1 = \xi_1 + \xi_2$ is   symmetric.

4. Let $a<0$ and $a\ne -1$.  In this case we construct   four  different types of distributions  $\mu_j$.

$4a$. Consider the subgroup  $H={\rm Ker}(I+\alpha)$. Then $H=\left\{\left(0, z\right): z\in\mathbb{T}\right\}= \mathbb{T}$. In this case we argue   as in case $3a$.

$4b$. Let $\gamma_j$  be Gaussian distributions on the group $\mathbb{R}$ with the characteristic functions  $\hat\gamma_j(s)=\exp\{-\sigma_js^2\}$, $j=1, 2$. Assume that
\begin{equation}\label{09_01_18_1}
\sigma_1+a\sigma_2=0.
\end{equation}
Consider the groups $\mathbb{R}$ and   $\mathbb{T}$ as subgroups of $X$ and put  $\mu_1=\gamma_1*m_\mathbb{T}$, $\mu_2=\gamma_2*\lambda$, where $\lambda$ is the distribution on the circle group $\mathbb{T}$ constructed in   case $3b$. Then the characteristic functions  $\hat\mu_j(s, n)$ satisfy Eq. (\ref{08_11_3}). Indeed,  taking into account that $\hat\mu_j(s, 0)=\hat\gamma_j(s)$, $j=1, 2$, (\ref{09_01_18_1}) implies that (\ref{08_11_3})  holds for $n_1=n_2=0$. If either $n_1\ne 0$ or $n_2\ne 0$, then   (\ref{08_01_18_2})  implies that both parts of Eq.  (\ref{08_11_3}) are equal to zero.

$4c$. Let $d\ge{b^2\over(a+1)^2}$ and $\sigma>0$. Consider the matrices
$$
A=\sigma \left(\begin{matrix}-a&-{ab\over a+1}\\ -{ab\over a+1}&\
d-{b^2\over a+1}\end{matrix}\right), \quad B=\sigma \left(\begin{matrix}1&{b\over a+1}\\ {b\over a+1}&\
d\end{matrix}\right).
$$
It is easy to check that the matrices $A$ and $B$ are  positive semidefinite  and the equality
\begin{equation}\label{10_11_1}
A+B\tilde\alpha=0
\end{equation}
holds. Denote by $\langle.,.\rangle$ the scalar product in $\mathbb{R}^2$. Consider Gaussian distributions  $\mu_1$ and $\mu_2$ on the group $X$  with the characteristic functions
$$
\hat\mu_1(s, n)=\exp\{-\langle A(s, n), (s, n)\rangle\}, \ \hat\mu_2(s, n)=\exp\{-\langle B(s, n), (s, n)\rangle\}, \ s\in\mathbb{R}, \ n\in\mathbb{Z}.
$$
It follows from (\ref{10_11_1}) that the characteristic functions $\hat\mu_j(s, n)$ satisfy Eq. (\ref{08_11_3}).

$4d$. Define the continuous monomorphism  ${\theta\colon\mathbb{R}\times\mathbb{Z}(2)\to X}$ by the formula
$$
\theta(t, k)=\left(t, (-1)^k e^{ibt\over {a+1}}\right), \quad t\in\mathbb{R}, \ k\in \mathbb{Z}(2).
$$
Put $S=\theta(\mathbb{R}\times\mathbb{Z}(2))$. Then $S$ is a close subgroup of   $X$, and $S$ is topologically isomorphic to   $\mathbb{R}\times\mathbb{Z}(2)$. It is easy to see that the restriction of  the automorphism $\alpha$ to the subgroup $S$ is a topological automorphism of the group   $S$. Let  $\lambda_j\in\Theta$, $j=1, 2$. Assume that   the characteristic functions  $\hat \lambda_j(s, l)$ are of the form
$$
\hat \lambda_j(s, l) = \begin{cases}\exp\{-\sigma_j s^2\},
&s\in \mathbb{R}, \ l=0,\\ \kappa_j\exp\{-\sigma'_j s^2\},
&s\in \mathbb{R}, \   l=1.
\\
\end{cases}
$$
Assume   that the equalities   $\sigma_1+a\sigma_2=0$ and $\sigma'_1+a\sigma'_2=0$ hold. Put $\mu_j=\theta(\lambda_j)$.
It is not difficult to verify that the characteristic functions  $\hat\mu_j(s, n)$ satisfy Eq. (\ref{08_11_3}).

 Let $\xi_1$ and
$\xi_2$ be independent random variables with values in the group
       $X$  and distributions
 $\mu_1$ and $\mu_2$, where $\mu_j$ constructed in cases $4b$--$4d$.   Since the characteristic functions  $\hat\mu_j(s, n)$ satisfy Eq. (\ref{08_11_3}), by Lemma  \ref{lem1}, the conditional distribution of the linear form
$L_2 = \xi_1 + \alpha\xi_2$
 given $L_1 = \xi_1 + \xi_2$ is symmetric.

Obviously, cases 1--4   exhaust   all possibilities for $a$ and $b$, when $\alpha=\left(\begin{matrix}a&b\\ 0& -1\end{matrix}\right)\in{\rm
Aut}(X)$. The examples of distributions $\mu_j$ constructed in case   3 and case 4  show that when either $a=-1$ and $b\ne 0$ or $a<0$ and $a\ne -1$  we can hardly expect
to get a reasonable description of distributions  which are characterized by the symmetry of the conditional distribution of the linear form   $L_2=\xi_1+\alpha\xi_2$ given $L_1=\xi_1+\xi_2$.

From what has been said above, when we considered cases 1--4, and Theorem \ref{th2}, we get the following assertion.
\begin{proposition}\label{pr2} Let $X=\mathbb{R}\times\mathbb{T}$. There is no   topological automorphism $\alpha$ of the group $X$ such that the following statements hold.

 $1$. If $\xi_1$ and $\xi_2$ are independent random variables with values in the group
       $X$  and distributions
 $\mu_1$ and $\mu_2$ with nonvanishing characteristic functions, then the symmetry of the conditional distribution of the linear form
$L_2 = \xi_1 + \alpha\xi_2$
 given $L_1 = \xi_1 + \xi_2$ implies that $\mu_j=\gamma_j*\rho_j$,
 where $\gamma_j$ are Gaussian distributions on $X$, and $\rho_j\in{\rm M}^1(G)$.

 $2$. There are independent random variables $\xi_1$ and
$\xi_2$  with values in the group
       $X$  and distributions  $\mu_j=\gamma_j*\rho_j$,
 where $\gamma_j$ are nondegenerate   Gaussian distributions on $X$  and $\rho_j\in{\rm M}^1(G)$, such that the conditional distribution of the linear form
$L_2 = \xi_1 + \alpha\xi_2$
 given $L_1 = \xi_1 + \xi_2$ is symmetric.
\end{proposition}

It is interesting to note that if $X={\mathbb
T}^2$ and $G$ is the subgroup of $X$ generated by all elements
of $X$ of order $2$,  then there is a topological automorphism   $\alpha$
 of the group $X$ such that statements  1 and   2 of Proposition \ref{pr2} are fulfilled  (see \cite{Fe4}).

\section{      Heyde's characterization theorem for the group $\Sigma_{{\text{\boldmath $a$}}}\times \mathbb{T}$}

Recall the definition of an $\text{\boldmath $a$}$-adic solenoid  $\Sigma_{\text{\boldmath $a$}}$. Let $\Delta_{{\text{\boldmath $a$}}}$ be
the group of ${\text{\boldmath $a$}}$-adic integers, where {{\text{\boldmath $a$}}}=$(a_0, a_1,\dots)$ is a  sequence of integers all greater than 1.   Consider the group
$\mathbb{R}\times\Delta_{{\text{\boldmath $a$}}}$, and denote by
 $B$ the subgroup of
$\mathbb{R}\times\Delta_{{\text{\boldmath $a$}}}$ of the form
$B=\{(n,n\mathbf{u}): n\in \mathbb{Z}\}$, where
$\mathbf{u}=(1, 0,\dots,0,\dots)$. The factor-group $\Sigma_{{\text{\boldmath $a$}}}=(\mathbb{R}\times\Delta_{{\text{\boldmath $a$}}})/B$ is called
  an ${\text{\boldmath $a$}}$-{adic solenoid}.  The group $\Sigma_{{\text{\boldmath $a$}}}$ is  compact, connected and has
dimension 1   (\!\!\cite[(10.12), (10.13),
(24.28)]{Hewitt-Ross}). The character group of the group   $\Sigma_{{\text{\boldmath $a$}}}$ is topologically isomorphic
to the discrete group of the form
$$H_{\text{\boldmath $a$}}= \left\{{m \over a_0a_1 \dots a_n} : \ n = 0, 1,\dots; \ m
\in {\mathbb{Z}} \right\}.
$$
We   identify $H_{\text{\boldmath $a$}}$ with the character group of the group $\Sigma_{{\text{\boldmath $a$}}}$.
It is convenient for us to consider $H_{\text{\boldmath $a$}}$  as a subset of $\mathbb{R}$.

Let $X$ be a   locally compact Abelian group, and  $n$ be an integer, $n\ne 0$. Denote by $f_n\colon X \to X$ an endomorphism of
the group $X$, defined by the formula  $f_nx=nx$, $x\in X$.
Each topological automorphism $a$ of the group
$\Sigma_\text{\boldmath $a$}$  can be written in the form
$a = f_p f_q^{-1}$ for some mutually prime $p$ и $q$,
where $f_p, f_q \in {\rm
Aut}(\Sigma_\text{\boldmath $a$})$. We   identify $a$
with the real number ${p\over q}$. Since $\tilde a =a$, we   also identify  $\tilde a$ with the real number ${p\over q}$.

Consider the group $X=\Sigma_{\text{\boldmath $a$}}\times\mathbb{T}$. The group $Y$  is topologically isomorphic to the group $H_{\text{\boldmath $a$}}\times\mathbb{Z}$.
  In order not to complicate the notation we assume that $Y=H_{\text{\boldmath $a$}}\times\mathbb{Z}$.
   Denote by $x=(g,z)$,
$g\in\Sigma_{\text{\boldmath $a$}}$, $z\in\mathbb{T}$, elements of the group $X$, and by $y=(r,n)$, $r\in H_{\text{\boldmath $a$}}$, $n\in\mathbb{Z}$, elements of the group $Y$. We   assume that the group $Y$ is embedded in a natural way into the group $\mathbb{R}\times \mathbb{Z}$. It is easy to verify that  each  automorphism $\alpha\in{\rm Aut}(X)$  is defined
by a matrix
$\left(\begin{matrix}a &b\\ 0&\pm 1\end{matrix}\right)$, where  $a\in
{\rm Aut}(\Sigma_\text{\boldmath $a$})$, $b\in H_{\text{\boldmath $a$}}$, and $\alpha$
is given by the formula
$$\alpha(g,z)=({a} g, (g,b)z^{\pm 1}), \quad
g\in\Sigma_\text{\boldmath $a$},  \ z\in\mathbb{T}.$$
Then the adjoint automorphism
$\tilde\alpha\in{\rm
Aut}(Y)$ is of the form
$$\tilde\alpha(r,n)=(a r+b n,\pm n),\quad r\in H_\text{\boldmath $a$}, \ n\in \mathbb{Z}.$$ We   identify the automorphisms ${\alpha}$ and ${\tilde\alpha}$
with the corresponding
matrix  $\left(\begin{matrix}a &b\\ 0&\pm
1\end{matrix}\right)$.

 It follows from  (\ref{f1}) and (\ref{f2}) that the characteristic function of a Gaussian distribution   $\gamma$ on   an  \text{\boldmath $a$}-adic solenoid
 $\Sigma_\text{\boldmath $a$}$ is of the form
$$
 \hat\gamma(r)=(g, r)\exp\{-\sigma r^2\}, \quad r\in H_{\text{\boldmath $a$}},
$$
where $g\in  \Sigma_\text{\boldmath $a$}$, $\sigma \ge 0$.

In this section, based on Theorem \ref{th2},   we will prove   a statement which can be viewed   as an analogue of   Heyde's theorem for the group  $X=\Sigma_{\text{\boldmath $a$}}\times\mathbb{T}$ for  two  independent random variables $\xi_1$ and $\xi_2$ with values in $X$ and linear forms
  $L_1=\xi_1+\xi_2$, $L_2=\xi_1+\alpha\xi_2$, where
  $\alpha\in{\rm
Aut}(X)$, and $\alpha$ is of the form $\alpha=\left(\begin{matrix}a&b\\ 0& 1\end{matrix}\right)$.
\begin{theorem}\label{th5}Let $X=\Sigma_{\text{\boldmath $a$}}\times\mathbb{T}$, where the
\text{\boldmath $a$}-adic solenoid  $\Sigma_\text{\boldmath $a$}$
contains no elements  of order $2$. Let $G$ be the subgroup of $X$ generated by the element of order $2$.    Let    $\alpha\in{\rm
Aut}(X)$, and let $\alpha=\left(\begin{matrix}a&b\\ 0& 1\end{matrix}\right)$, where ${\rm Ker}(I+a)=\{0\}$.
Let $\xi_1$ and $\xi_2$ be independent random variables with values in the group
       $X$  and distributions
 $\mu_1$ and $\mu_2$ with nonvanishing characteristic functions.
Assume that the conditional distribution of the linear form
$L_2=\xi_1+\alpha\xi_2$
 given $L_1 = \xi_1 + \xi_2$ is symmetric. If  $a
< 0$, then there is a continuous monomorphism
$\phi\colon\mathbb{R}\times\mathbb{Z}(2)\to X$ such that
the distributions $\mu_{j}$ are represented in the form $\mu_j=\phi(M_j)*E_{x_j}$,
  where $x_j\in X$, and $M_j$ are distributions on the group  $\mathbb{R}\times\mathbb{Z}(2)$   belonging to the class $\Theta$.
  If $a>0$,  then some shifts of the distributions  $\mu_j$
are supported in $G$.
\end{theorem}
To prove Theorem \ref{th5} we need the following lemma, which deals with arbitrary locally compact Abelian groups.
\begin{lemma} [\!\!{\protect\cite{FeTVP1}}] \label{newlem6}
    Let  $X$ be a locally compact Abelian group containing
  no elements of order $2$.  Let $\alpha$ be a topological automorphism of the group
$X$ satisfying the condition   ${\rm Ker}(I+\alpha)=\{0\}$.
Let  $\xi_1$ and  $\xi_2$ be independent random variables with values in
$X$   and distributions   $\mu_1$ and $\mu_2$ with nonvanishing characteristic functions. Then the symmetry of the conditional distribution of the linear form
$L_2 = \xi_1 + \alpha\xi_2$ given  $L_1 = \xi_1 +
\xi_2$ implies that   $\mu_1$ and $\mu_2$ are Gaussian distributions.
\end{lemma}

{\it Proof of Theorem $\ref{th5}$} \
By Lemma \ref{lem1}, the characteristic functions  $\hat\mu_j(y)$ satisfy Eq. (\ref{42}). Put
$\nu_j = \mu_j
* \bar \mu_j$. Then
 $\hat \nu_j(y) = |\hat \mu_j(y)|^2 > 0,$   $y \in Y$, $j=1, 2$.
 Obviously, the characteristic functions  $\hat \nu_j(y)$
also satisfy Eq. (\ref{42}) which takes the form
\begin{equation}\label{21_09_1}
\hat\nu_1(r_1+r_2, n_1+n_2)\hat\nu_2(r_1+a r_2+bn_2, n_1+n_2)$$$$=
\hat\nu_1(r_1-r_2, n_1-n_2)\hat\nu_2(r_1-a r_2-bn_2, n_1-n_2), \quad r_1, r_2\in H_{\text{\boldmath $a$}},
 \ n_1, n_2\in \mathbb{Z}.
\end{equation}
Substituting  $n_1=n_2=0$ in
(\ref{21_09_1})  we get
\begin{equation}\label{21_09_2}
\hat\nu_1(r_1+r_2, 0)\hat\nu_2(r_1+a r_2,0)=
\hat\nu_1(r_1-r_2, 0)\hat\nu_2(r_1-a r_2, 0), \quad r_1, r_2\in H_{\text{\boldmath $a$}}.
\end{equation}
Applying Lemmas  \ref{lem1} and   \ref{newlem6}  to the \text{\boldmath $a$}-adic solenoid $\Sigma_\text{\boldmath $a$}$, we obtain from (\ref{21_09_2}) that
\begin{equation}\label{21_09_3}
\hat\nu_j(r, 0)=e^{-\sigma_j r^2}, \quad r\in H_{\text{\boldmath $a$}},
\end{equation}
where $\sigma_j\ge 0$,   $j=1, 2$.

Note that the inequality
$$
|\hat\nu(u)-\hat\nu(v)|^2\le 2(1-{\rm Re}\
\hat\nu(u-v)), \quad u, v\in Y,
$$
holds for each characteristic function $\hat\nu(y)$ on the group  $Y$. Taking this into account, it follows from (\ref{21_09_3})   that   the characteristic functions   $\hat\nu_j(r, n)$ are uniformly continuous on   $Y$ in the topology induced on $Y$ by the
topology of the group
$\mathbb{R}\times \mathbb{Z}$. This implies that   the characteristic functions $\hat\nu_j(r, n)$ can
be extended by continuity to some continuous positive definite functions $g_j(s, n)$ on the group $\mathbb{R}\times \mathbb{Z}$. By the Bochner theorem, there are distributions $\lambda_j$  on the group
 $\mathbb{R}\times \mathbb{T}$ such that $\hat \lambda_j(s, n)=g_j(s, n)$, $s\in\mathbb{R}$, $n\in\mathbb{Z}$, $j=1, 2$.
Denote by $\iota$ the natural embedding of the group  $Y=H_{\text{\boldmath $a$}}\times \mathbb{Z}$ into the group $\mathbb{R}\times \mathbb{Z}$. Put $\upsilon =\tilde \iota$. Then
$\upsilon \colon \mathbb{R}\times \mathbb{T}\to \Sigma_{\text{\boldmath $a$}}\times \mathbb{T}$ is a continuous homomorphism, and $\nu_j=\upsilon (\lambda_j)$, $j=1, 2$.
Since the image
$\iota(Y)$ is dense in the group $\mathbb{R}\times \mathbb{Z}$,
$\upsilon $ is a monomorphism.

It follows from (\ref{21_09_1}) that the characteristic functions  $\hat \lambda_j(s, n)$
satisfy the equation
\begin{equation}\label{23_1}
\hat \lambda_1(s_1+s_2, n_1+n_2)\hat \lambda_2(s_1+a s_2+bn_2, n_1+n_2)$$$$=
\hat \lambda_1(s_1-s_2, n_1-n_2)\hat \lambda_2(s_1-a s_2-bn_2, n_1-n_2), \quad s_1, s_2\in\mathbb{R}, \ n_1, n_2\in\mathbb{Z}.
\end{equation}
We retain the notation $G$ for the subgroup of $\mathbb{R}\times \mathbb{T}$
generated by the element of order $2$ and the notation $\alpha$ for the topological automorphism of the group
$\mathbb{R}\times \mathbb{T}$ which is defined by the matrix $\left(\begin{matrix}a&b\\ 0& 1\end{matrix}\right)$.  Let
$\eta_1$ and  $\eta_2$ be independent random variables with values in the group
 $\mathbb{R}\times \mathbb{T}$ and distributions  $\lambda_1$ and $\lambda_2$.
By Lemma \ref{lem1},   (\ref{23_1}) implies that the conditional distribution of the linear form
 $P_2 = \eta_1 + \alpha\eta_2$   given
 $P_1 = \eta_1 + \eta_2$ is symmetric.

 Let $a<0$. It follows from ${\rm Ker}(I+a)=\{0\}$ that $a\ne -1$. Since the characteristic functions $\hat \nu_j(y)$ do not vanish, by Theorem \ref{th2}, $\lambda_j=\tau(M_j)$, where a continuous monomorphism     $\tau\colon\mathbb{R}\times\mathbb{Z}(2)\to \mathbb{R}\times\mathbb{T}$ is defined by  formula  (\ref{22_1}), and   $M_j\in \Theta$.
 Put $\phi=\upsilon \tau$. Then $\phi\colon\mathbb{R}\times\mathbb{Z}(2)\to \Sigma_{\text{\boldmath $a$}}\times\mathbb{T}$ is a continuous monomorphism.
 It follows from the Suslin theorem that images of Borel sets under the mapping
$\phi$ are also Borel sets, so that $\phi$ induces an isomorphism of the convolution
semigroups ${\rm M}^1(\mathbb{R}\times \mathbb{Z}(2))$ and ${\rm M}^1(\phi(\mathbb{R}\times \mathbb{Z}(2))$.
We have $\nu_j=\phi(M_j)$, $j=1, 2$.
Applying Lemmas  \ref{lem11}  and   \ref{lem8}, we get that some shifts of the distributions
 $\mu_j$ are   images under the action of $\phi$ of distributions on  the group $\mathbb{R}\times\mathbb{Z}(2)$   belonging to the class $\Theta$. Thus, in the case when   $a<0$, the theorem is proved.

Let $a>0$. Substituting (\ref{21_09_3}) into  (\ref{21_09_2})  we get
$\sigma_1+a\sigma_2=0$. Hence  $\sigma_1=\sigma_2=0$,
and (\ref{21_09_3}) implies that $\hat\nu_j(r, 0)=1$ for all $r\in H_{\text{\boldmath $a$}}$, $j=1, 2$.
Applying Lemma \ref{lem2},  we get that $\sigma(\nu_j)\subset A(X, H_{\text{\boldmath $a$}})=\mathbb{T}$.
  Let
$\zeta_j$ be independent random variables with values in the group
 $X$ and distributions $\nu_j$. By Lemma \ref{lem1}, it follows from (\ref{21_09_1}) that the conditional distribution of the linear form $Q_2 = \zeta_1 + \alpha\zeta_2$ given
  $Q_1 = \zeta_1 + \zeta_2$ is symmetric.
 Note that the restriction of the automorphism $\alpha$ to the subgroup $\mathbb{T}$
 is the identity automorphism. Thus, if we consider  $\zeta_j$ as
 independent random variables with values in the circle group
 $\mathbb{T}$, then
 the conditional distribution of the linear form  $Q_2 = \zeta_1 + \zeta_2$ given
  $Q_1 = \zeta_1 + \zeta_2$ is symmetric. By Lemmas  \ref{lem1} and    \ref{lem2}, this   implies that
 $\sigma(\nu_j)\subset G$. By Lemma \ref{lem11}, some shifts of the distributions   $\mu_j$
are supported in $G$.    $\blacksquare$
\medskip

Let $X=\Sigma_{\text{\boldmath $a$}}\times\mathbb{T}$, where the
\text{\boldmath $a$}-adic solenoid  $\Sigma_\text{\boldmath $a$}$
  does not contain an element  of order $2$. We note that since   $\Sigma_{\text{\boldmath $a$}}$ is a connected compact Abelian group,   multiplication by 2 is an epimorphism of the  group    $\Sigma_{\text{\boldmath $a$}}$. It follows from this that if  the group $\Sigma_\text{\boldmath $a$}$
  does not contain an element  of order $2$, then $\Sigma_\text{\boldmath $a$}$ and $H_\text{\boldmath $a$}$ are groups with unique division by 2. Let  $\alpha\in{\rm
Aut}(X)$, and  let $\alpha=\left(\begin{matrix}a&b\\ 0& 1\end{matrix}\right)$.  Omit in Theorem \ref{th5}   the condition ${\rm Ker}(I+a)=\{0\}$, and consider the case when
 ${\rm Ker}(I+a)=\Sigma_{\text{\boldmath $a$}}$. The following statement holds.
 \begin{proposition}\label{pr4}
 Let $X=\Sigma_{\text{\boldmath $a$}}\times\mathbb{T}$, where the
\text{\boldmath $a$}-adic solenoid  $\Sigma_\text{\boldmath $a$}$
 does not   contain an element   of order $2$. Let $G$ be the subgroup of $X$ generated by the element of order $2$.    Let    $\alpha=\left(\begin{matrix}-1&b\\ 0& 1\end{matrix}\right)\in{\rm
Aut}(X)$. Let $K$ be a subgroup of $X$ of the form
$K=\left\{\left(g, \pm {\left(-g, {b\over 2}\right)}\right): g\in \Sigma_{\text{\boldmath $a$}}\right\}$. Let $\xi_1$ and $\xi_2$ be independent random variables with values in the group
       $X$  and distributions
 $\mu_1$ and $\mu_2$.
Assume that the conditional distribution of the linear form
$L_2=\xi_1+\alpha\xi_2$
 given $L_1 = \xi_1 + \xi_2$ is symmetric.
 Then the distributions $\mu_j$  are of the form $\mu_j=\lambda_j*E_{x_j},$
where   $\lambda_j\in{\rm M}^1(K)$, $x_j\in X$.  Moreover,
either $\lambda_2=\lambda_1*\omega_1$ or $\lambda_1=\lambda_2*\omega_2$, where $\omega_j\in {\rm M}^1(G)$.
\end{proposition}
{\it Proof} \ We observe that $G=\{(0, \pm 1)\}\subset K$. Put
 $L={\rm Ker}(I-\tilde\alpha)$.   Obviously,  $L=\{(r, n):2r-bn=0\}=\left\{\left({bn\over 2}, n\right): n\in \mathbb{Z}\right\}$, and hence $L^{(2)}=\{(bn, 2n): n\in \mathbb{Z}\}$. It follows from this that $K=A\left(X, L^{(2)}\right)$. Arguing as in the proof of Lemma  \ref{lem23.1}, we obtain that   $\hat\mu_j(y)=(p_j, y)$, $y\in L^{(2)}$, $p_j\in X$, $j=1, 2$. Moreover $p_1+p_2\in K$. Let $p_j=(g_j, z_j)$, $j=1, 2$. Put $\tilde x=\left(g_1,  \left(-g_1, {b\over 2}\right)\right)$,  $x_1=p_1-\tilde x$, $x_2=-p_1+\tilde x$, $\lambda_j=\mu_j*E_{-x_j}$, $j=1, 2$. Denote by $\eta_j$ independent random variables with values in the group
       $X$  and distributions
$\lambda_j$. As in the proof of Lemma  \ref{lem23.1}, we verify that the conditional distribution of the linear form    $N_2=\eta_1+\alpha\eta_2$ given $N_1=\eta_1+\eta_2$   is symmetric, and $\sigma(\lambda_j)\subset K$, $j=1, 2$. It is easy to see that the  automorphism   $\alpha$  acts on $K$ as $-I$.  Thus, if we consider $\eta_j$ as independent random variables with values in the group $K$, then the conditional distribution of the linear form    $N_2=\eta_1-\eta_2$ given $N_1=\eta_1+\eta_2$   is symmetric. Applying Proposition \ref{pr1} to the group $K$ and independent random variables $\eta_j$, we obtain that either $\lambda_2=\lambda_1*\omega_1$ or $\lambda_1=\lambda_2*\omega_2$, where $\omega_j\in {\rm M}^1(G)$, $j=1, 2$. The proposition   follows from  this. $\blacksquare$
\medskip

Theorems \ref{th2} and   \ref{th5} imply   the following statement.
\begin{proposition}\label{pr3}  Let either $X=\mathbb{R}\times\mathbb{T}$  or $X=\Sigma_{\text{\boldmath $a$}}\times\mathbb{T}$, where the
\text{\boldmath $a$}-adic solenoid  $\Sigma_\text{\boldmath $a$}$
does not contain  an element   of order $2$.  Let $\alpha\in{\rm
Aut}(X)$. Denote by $G$ the subgroup of $X$ generated by the element of order $2$. Assume that
\begin{equation}\label{25_01_1}
{\rm Ker}(I+\alpha)=G.
\end{equation}
Let $\xi_1$ and $\xi_2$ be independent random variables with values in the group
       $X$  and distributions
  $\mu_1$ and $\mu_2$ with nonvanishing characteristic functions.
If the conditional distribution of the linear form
$L_2=\xi_1+\alpha\xi_2$
 given $L_1 = \xi_1 + \xi_2$ is symmetric, then there is a continuous monomorphism
$\pi\colon\mathbb{R}\times\mathbb{Z}(2)\to X$ such that
the distributions $\mu_{j}$ are represented in the form $\mu_j=\pi(M_j)*E_{x_j}$,
  where  $M_j$ are distributions on the group  $\mathbb{R}\times\mathbb{Z}(2)$   belonging to the class $\Theta$, and $x_j\in X$.
\end{proposition}
{\it Proof} \  It follows from (\ref{25_01_1}) that $\alpha=\left(\begin{matrix}a&b\\ 0& 1\end{matrix}\right)$. Really, suppose $\alpha=\left(\begin{matrix}a&b\\ 0& -1\end{matrix}\right)$. Then ${\rm Ker}(I+\alpha)\supset\mathbb{T}$, contrary to (\ref{25_01_1}).
Let the continuous monomorphisms   $\tau$  and  $\phi$ be the same as in  Theorems  \ref{th2} and   \ref{th5} respectively.

Let $X=\mathbb{R}\times\mathbb{T}$. Then (\ref{25_01_1}) holds if and only if $a\ne-1$. Indeed, if $a=-1$, then ${\rm Ker}(I+\alpha)=F\times G$, contrary to (\ref{25_01_1}). The converse is obvious. Thus, $a\ne -1$. Taking into account that  the characteristic functions $\hat\mu_j(y)$ do not vanish,   in Theorem \ref{th2} only cases  $1b$ and   $3$ are possible, and we can put $\pi=\tau$.

Let $X=\Sigma_{\text{\boldmath $a$}}\times\mathbb{T}$. Then (\ref{25_01_1}) holds if and only if ${\rm Ker}(I+a)=\{0\}$. Indeed,
if ${\rm Ker}(I+a)\ne\{0\}$, then   $G$ is a proper   subgroup of ${\rm Ker}(I+\alpha)$, contrary to (\ref{25_01_1}). The converse is obvious. Thus,   ${\rm Ker}(I+a)=\{0\}$, we are in the conditions of Theorem \ref{th5} and can put $\pi=\phi$. $\blacksquare$

It is interesting to note that  if $X=\Sigma_{\text{\boldmath $a$}}$, where the $\text{\boldmath $a$}$-adic solenoid  $\Sigma_{\text{\boldmath $a$}}$ contains  an element of order 2, then Proposition \ref{pr3} holds for the group $X$ (\!\!\cite{F_solenoid}).

\vskip 1 cm

\noindent B. Verkin Institute for Low Temperature Physics and Engineering\\
of the National Academy of Sciences of Ukraine\\
47, Nauky Ave, Kharkiv, 61103, Ukraine

\noindent e-mail: feldman@ilt.kharkov.ua

\end{document}